\documentclass{svjour3}
\usepackage{amssymb}
\usepackage{url}

\usepackage{graphicx}
\usepackage{amstext,amsmath,amscd}
\usepackage{color}
\usepackage[curve]{xypic}
\usepackage{psfrag}

\definecolor{DkBlue}{rgb}{0,0,.5}

\newcommand{\D} {\mathbb D}
\newcommand{\I} {\mathcal I}
\newcommand{\M} {\mathbb M}
\newcommand{\Mi}[1]{{\mathbb M}_{#1}^{-1}}
\newcommand{\p} {\partial}
\newcommand{\R} {\mathbb R}
\newcommand{\Rmap}{\mathcal{P}}

\newcommand{\bp}{\textbf{p}}
\newcommand{\bq}{\textbf{q}}
\newcommand{\bv}{\textbf{v}}
\newcommand{\bx}{\textbf{x}}
\newcommand{\by}{\textbf{y}}

\newcommand{\ds}{\displaystyle}
\newcommand{\Div}{\text{div }}
\newcommand{\Curl}{\text{curl }}
\newcommand{\pc}[1]{\textsc{#1}}
\newcommand{\dc}[1]{{\overline{\textsc{#1}}}}
\newcommand{\cs}{\mathcal{C}}
\newcommand{\dcs}{\overline{\mathcal{C}}}
\newcommand{\ol}[1]{\overline{#1}}

\newcommand{\Tau}{\mathcal T}

\newcommand{\cancel}[1]{}
\newcommand{\raw}{\rightarrow}

\newcommand{\whitC}[3]{\lambda_{#1}\nabla \lambda_{#2}\times\nabla \lambda_{#3}}

\newcommand{\vsn}[1]{\left|#1\right|}

\newcommand{\wabwab}{\alpha}
\newcommand{\wabwac}{\beta}
\newcommand{\wacwac}{\gamma}
\newcommand{\wacwbc}{\delta}

\newcommand{\wabwabn}{{\frac{12P^2+1}{24P}}}
\newcommand{\wabwacn}{{\frac{4P^2-1}{48P}}}
\newcommand{\wacwacn}{{\frac{12 P^2 + 20\sqrt{3} P + 21}{144 P}}}
\newcommand{\wacwbcn}{{\frac{4P^2-5}{48P}}}

\newcommand{\dwabwab}{\vartheta}
\newcommand{\dwabwac}{\zeta}
\newcommand{\dwacwac}{\theta}
\newcommand{\dwacwad}{\kappa}
\newcommand{\dwacwbc}{\xi}

\newcommand{\dwabwabn}{\left(\eta_{\star\sigma^1_{12}},\eta_{\star\sigma^1_{12}}\right)}
\newcommand{\dwabwacn}{\left(\eta_{\star\sigma^1_{12}},\eta_{\star\sigma^1_{13}}\right)}
\newcommand{\dwacwacn}{\left(\eta_{\star\sigma^1_{13}},\eta_{\star\sigma^1_{13}}\right)}
\newcommand{\dwacwadn}{\left(\eta_{\star\sigma^1_{13}},\eta_{\star\sigma^1_{14}}\right)}
\newcommand{\dwacwbcn}{\left(\eta_{\star\sigma^1_{13}},\eta_{\star\sigma^1_{23}}\right)}

\newcommand{\mdiagent}{\varrho}
\newcommand{\mdiagentn}{\frac 1{4P^4}+\frac{P}{\sqrt{3+12P^2}}}

\newcommand{\cond}[1]{\text{cond}\left(#1\right)}

\newcommand{\whit}{\mathcal{W}}
\newcommand{\dwhit}{\ol{\mathcal{W}}}

\newcommand{\Mdiag}[1]{{\M^{Diag}_{#1}}}
\newcommand{\Mwhit}[1]{{\M^{Whit}_{#1}}}

\newcommand{\Mduali}[1]{{\left(\M^{Dual}_{#1}\right)^{-1}}}

\newcommand{\bc}{\textbf{c}}

\newcommand{\hcurl}{H(\textnormal{curl})}

\newtheorem{defn}{Definition}
\newtheorem{cor}{Corollary}

\begin{document}

\thispagestyle{empty}
\begin{center}
{\LARGE\textbf{Dual Formulations of Mixed Finite Element Methods with Applications}} \\
\textsc{Andrew Gillette\footnote{Department of Mathematics, University of Texas at Austin, \url{agillette@math.utexas.edu}}, Chandrajit Bajaj\footnote{Department of Computer Science, Institute for Computational Engineering and Sciences, University of Texas at Austin, \url{bajaj@cs.utexas.edu}}} \\
\today
\end{center}

\begin{abstract}
Mixed finite element methods solve a PDE using two or more variables.  The theory of Discrete Exterior Calculus explains why the degrees of freedom associated to the different variables should be stored on both primal and dual domain meshes with a discrete Hodge star used to transfer information between the meshes.  We show through analysis and examples that the choice of discrete Hodge star is essential to the numerical stability of the method.  Additionally, we define interpolation functions and discrete Hodge stars on dual meshes which can be used to create previously unconsidered mixed methods.  Examples from magnetostatics and Darcy flow are examined in detail.
\end{abstract}

%Keywords:
%Discrete exterior calculus \sep Finite element method \sep Partial differential equations \sep Whitney forms \sep Hodge star

\section{Introduction}

The theory of Discrete Exterior Calculus (DEC) has provided a novel viewpoint for analyzing linear systems derived from finite element theory.  We highlight three important conclusions of this theory:
\begin{enumerate}
\item Variables in a PDE should be discretized as degree of freedom arrays (``cochains'') over a primal simplicial mesh or its dual mesh.
\item A discrete Hodge star is used to transfer information between primal and dual meshes.
\item Whitney elements provide stable finite elements for the primal mesh.
\end{enumerate}
Most numerical methods for PDEs over unstructured tetrahedral meshes discretize variables as cochains over the primal mesh and build up linear systems from there.  In this paper, we look at the alternative approach of discretizing variables over the dual mesh and design dual formulations of the linear systems based on DEC theory.  This approach is especially valuable in the context of mixed finite element systems as they employ all the key ingredients of DEC theory: both primal and dual cochains, a discrete Hodge star, and, typically, Whitney elements.

Before turning to mixed systems, however, we look at a simpler example from electromagnetics illustrating the relevance and benefit of our technique.  The example is inspired by He and Teixeira~\cite{HT2006}.  Using a Discrete Exterior Calculus analysis of Maxwell's equations, one can derive a second order vector wave equation
\begin{equation}
\label{eq:prwave}
\D_1^T\M_2\D_1\pc e = \omega^2 \M_1\pc e,
\end{equation}
where $\pc e$ is the electric field intensity, discretized as a cochain on the primal mesh, $\omega$ is a coefficient, $\D_1$ is a rectangular incidence matrix having entries of $0$ and $\pm 1$ only, and $\M_k$ is a discrete Hodge star operator.

The dual formulation of this physical phenomenon is an equation for the magnetic field intensity $\dc h$, discretized as a cochain on the dual mesh:
\begin{equation}
\label{eq:duwave}
\D_1\Mi 1\D_1^T\dc h = \omega^2 \Mi 2\dc h.
\end{equation}
Both systems (\ref{eq:prwave}) and (\ref{eq:duwave}) are computationally tractable if $\M_k$ is a diagonal matrix which, by DEC theory, can be achieved when the primal and dual meshes are orthogonal.  If orthogonality is not guaranteed, as is the case with barycentric dual meshes, $\M_k$ is defined using Whitney elements and results in a sparse matrix.  As a consequence, system (\ref{eq:duwave}) then involves possibly full rank matrices and is thus significantly more computationally expensive to solve.  He and Teixeira~\cite{HT2006} reduce the rank of the $\Mi k$ matrices by using a topological thresholding technique which requires an input parameter.

Our approach skirts the problem of full rank inverses by introducing a novel definition of the $\Mi k$ matrices free of parameters and guaranteed to produce a sparse matrix.  The outline of the paper and summary of its contributions are as follows:

\begin{itemize}

\item In Section \ref{sec:pwnot}, we briefly discuss prior work and fix relevant notation.

\item In Section \ref{sec:dualwhit}, we use the Sibson coordinate functions to construct dual Whitney-like functions which define a novel sparse inverse discrete Hodge star $(\M_k^{Dual})^{-1}$.  We show how the choice of discrete Hodge star requires certain geometric quality conditions of the primal and dual mesh elements.  A specific example is given showing how our dual formulation of the problem can result in a better conditioned linear system than the primal formulations.

\item In Section \ref{sec:res}, we examine how our methodology applies to generic PDE problems as well as to some specific applications employing mixed finite element methods.  We cast each into our common notational framework and show how to formulate equivalent dual formulations of the problem from a DEC-based analysis.  The specific advantages of these dual formulations are analyzed, including an ability to compare and contrast calculations on a primal mesh with the analogous calculations on the dual mesh.
\end{itemize}

\section{Prior Work and Notation}
\label{sec:pwnot}

Our work is inspired primarily by the emergent theory of Discrete Exterior Calculus (DEC).  DEC is an attempt to create from scratch a discrete theory of differential geometry and topology whose definitions and theorems mimic their continuous counterparts~\cite{H2003,DHLM2005}.  A central conclusion of the theory is that degrees of freedom for finite elements should be assigned to mesh vertices, edges, faces or interiors according to the dimensionality of the variable being modeled.  If these degrees of freedom have a natural geometric duality, as occurs for example between electric and magnetic fields, two meshes of the domain are necessary - a primal and dual mesh \cite{H2001b}.  This has given rise to DEC-based methods for solving problems of Darcy flow  \cite{HNC2008}, electromagnetism \cite{HT2006} and elasticity \cite{Y2008}, among others.  As we will show, the `bottom-up' approach of DEC clearly suggests alternative discretization methods less evident from such 'top-down' theories as finite element exterior calculus~\cite{AFW2010}.

\begin{figure}
\centering
\psfrag{sigma0}{$\sigma^0$}
\psfrag{starsigma0}{$\star\sigma^0$}
\psfrag{sigma1}{$\sigma^1$}
\psfrag{starsigma1}{$\star\sigma^1$}
\psfrag{sigma2}{$\sigma^2$}
\psfrag{starsigma2}{$\star\sigma^2$}
\includegraphics[width=.6\linewidth]{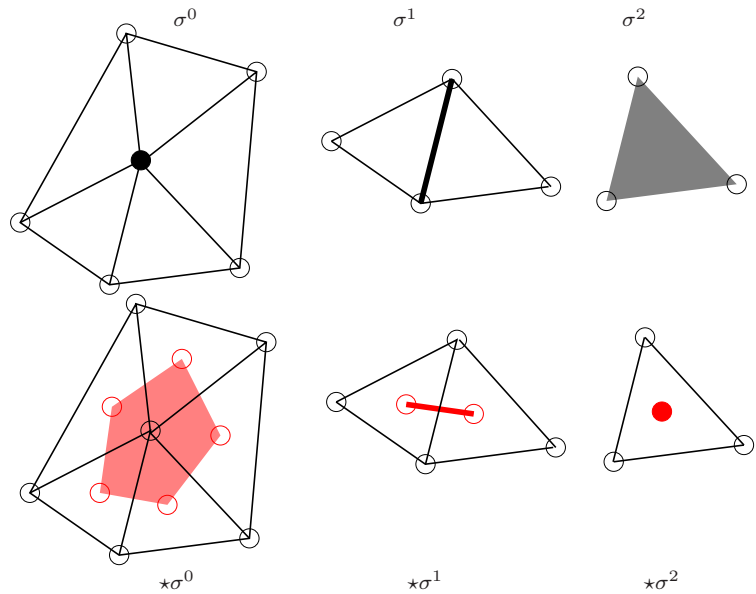}
\caption{Primal simplices are shown in black in the top row: $\sigma^0$ is a vertex, $\sigma^1$ is an edge, and $\sigma^2$ is a face.  Their corresponding dual cells for $n=2$ are shown in red on bottom: $\star\sigma^2$ is the barycenter of $\sigma^2$, $\star\sigma^1$ is an edge between barycenters, and $\star\sigma^0$ is a planar polygon with barycenters as vertices.  In three dimensions ($n=3$), primal vertices have dual polytopes, primal edges have dual polygonal facets, primal faces have dual edges, and primal volumes have dual vertices.}
\label{fig:2dprimaldual}
\end{figure}

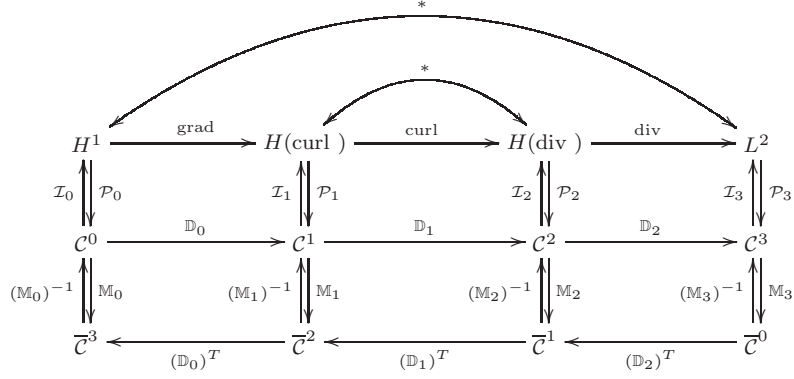
\begin{figure*}
\[
\xymatrix{
H^1 \ar[rr]^{\text{grad}} \ar@<.4ex>[d]^{\Rmap_0} \ar@{<->} @/^4pc/[rrrrrr]^\ast  && H(\Curl) \ar[rr]^{\Curl} \ar@<.4ex>[d]^{\Rmap_1}  \ar@{<->} @/^2pc/[rr]^\ast  && H(\Div)  \ar[rr]^{\Div} \ar@<.4ex>[d]^{\Rmap_2} && L^2 \ar@<.4ex>[d]^{\Rmap_3} \\
{\cs^0}  \ar[rr]^{\D_0} \ar@<.4ex>[u]^{\I_0} \ar@<.4ex>[d]^{\M_0} && {\cs^1}  \ar[rr]^{\D_1} \ar@<.4ex>[u]^{\I_1} \ar@<.4ex>[d]^{\M_1} && {\cs^2}   \ar[rr]^{\D_2} \ar@<.4ex>[u]^{\I_2} \ar@<.4ex>[d]^{\M_2} &&  {\cs^3} \ar@<.4ex>[u]^{\I_3} \ar@<.4ex>[d]^{\M_3} \\
{\dcs^3}  \ar@<.4ex>[u]^{(\M_0)^{-1}} && {\dcs^2}  \ar[ll]^{(\D_0)^T} \ar@<.4ex>[u]^{(\M_1)^{-1}}  && {\dcs^1} \ar[ll]^{(\D_1)^T}   \ar@<.4ex>[u]^{(\M_2)^{-1}} &&  {\dcs^0} \ar[ll]^{(\D_2)^T}  \ar@<.4ex>[u]^{(\M_3)^{-1}}
}
\]
\caption{The combined DEC and deRham diagram for a contractible domain in $\R^3$.  The top row shows the $L^2$ deRham diagram with continuous Hodge star maps between function spaces. The middle and bottom rows show primal and dual cochain spaces, respectively, along with the discrete exterior derivative and discrete Hodge star maps. The $\I$ and $\Rmap$ maps are interpolation (Whitney) and projection (deRham) maps. }
\label{fig:derham}
\end{figure*}

The main notational aspects of DEC are encapsulated by Figures \ref{fig:2dprimaldual} and \ref{fig:derham}.   Figure~\ref{fig:2dprimaldual} shows our notation for domain elements, i.e. primal $k$-simplices $\sigma^k$ and their geometric dual $n-k$-cells $\star\sigma^{n-k}$ where $n$ is the dimension of the domain.  The dual domain mesh is defined by taking the circumcenters or barycenters of $n$-simplices and connecting them based on simplex adjacency in the usual manner.  The measure of $\sigma^k$ (respectively $\star\sigma^{n-k}$) is denoted $|\sigma^k|$ (respectively $|\star\sigma^{n-k}|$), meaning length for $k=1$, area for $k=2$, and volume for $k=3$, with the convention that $|\sigma^0|=|\star\sigma^n|=1$.

Figure~\ref{fig:derham} shows the various continuous and discrete spaces relevant to DEC theory for $n=3$ and the operators between them.  The vector space of $k$-cochains, i.e. linear mappings from $k$-simplices to $\R$, is denoted $\cs^k$.  The vector space of dual $k$-cochains, i.e. linear mappings from $k$-cells of the dual mesh to $\R$, is denoted $\dcs^k$.  The $\D_k$ matrix is the transpose of the $(k+1)$st boundary operator, i.e. it encodes element adjacency and orientation information with entries $\pm 1$.

The interpolation map $\I_k$ converts a $k$-cochain into a piecewise-defined $k$-form whose global continuity in a distributional sense is indicated by Figure~\ref{fig:derham} (e.g. $\I_1\pc w\in\hcurl$).   Define $\I_k$ by
\begin{equation}
\label{eqn:ikdef}
\I_k(\pc w):=\sum_{{\sigma^k}\in\cs_k}\pc w({\sigma^k})\whit_{\sigma^k}.
\end{equation}
where $\whit_{\sigma^k}$ is the Whitney function associated to simplex $\sigma^k$.  These functions are described in Appendix A.  The Whitney functions were first described in~\cite{W1957} and later recognized by Bossavit \cite{B1988b} and others as the correct generalization of edge and face elements needed for DEC theory.  An extensive treatment of all of these spaces, functions, and operators is given in~\cite{G2011}.

We now discuss the Hodge star $\ast$ and its discretization as a square matrix $\M$ or $\M^{-1}$.  As shown in Figure~\ref{fig:derham}, the continuous Hodge star $\ast$ maps between forms of complementary and orthogonal dimensions, i.e. $\ast:\Lambda^k\raw\Lambda^{n-k}$.  For domains in $\R^3$ as considered here, $\ast$ is defined by the equations
\[\ast 1=dxdydz,\quad\ast dx=dydz,\quad\ast dy=-dxdz,\quad\ast dz=dxdy,\quad \ast\ast=1.\]
For a more general definition of $\ast$, see~\cite{AMR1988}.

A discrete Hodge star $\M$ maps not only between cochains of complementary dimensions but also between primal and dual meshes~\cite{H2001b}.  In this paper, we focus on the two definitions of a discrete Hodge star most relevant to DEC theory.  The first is the diagonal discrete Hodge star defined by
\begin{equation}
\label{eqn:hsdiag}
(\M_k^{Diag})_{ij} := \frac{|\star\sigma^k_i|}{|\sigma^k_i|}\delta_{ij}.
\end{equation}
The definition of $\M_k^{Diag}$ fits nicely into DEC theory when the dual mesh is defined by taking circumcenters of the primal simplices, thus producing orthogonal meshes~\cite{DHLM2005}.  In practice, however, it is often desirable to use barycenters to define the dual mesh as this guarantees that $\sigma^k$ will intersect $\star\sigma^k$ in the ambient space.  A correction factor for this change is given by Auchmann and Kurz~\cite{AK2006}.

The more widely used approach for barycentric dual meshes employs Whitney interpolants in the definition of the discrete Hodge star:
\begin{equation}
\label{eqn:hswhit}
(\M_k^{Whit})_{ij} := \left(\whit_{\sigma^k_i},\whit_{\sigma^k_j}\right) = \int_K \whit_{\sigma^k_i}\cdot\whit_{\sigma^k_j}
\end{equation}
The inner product here is the standard integration of scalar or vector valued functions over the domain.  Dodziuk~\cite{D1976} originally proposed the definition of $\M_k^{Whit}$ but it has been called the Galerkin Hodge~\cite{B2000part5} for its relation to finite element methods.  Bell~\cite{wB2008} has implemented linear solvers in a DEC context using $\M_k^{Whit}$ for various $k$.

Many other discrete Hodge stars appear in the literature, including the combinatorial discrete Hodge star of Wardetzsky and Wilson~\cite{W2006,W2007} and the metrized chain Hodge star of DiCarlo et al.~\cite{DMPS2009}.  To our knowledge, no authors have defined a discrete Hodge star using dual interpolatory functions as we propose in this work.

\section{Dual Whitney Interpolants and Dual Discrete Hodge Stars}
\label{sec:dualwhit}
It is evident from the DEC-deRham diagram in Figure~\ref{fig:derham} that the direct interpolation of degrees of freedom on a dual mesh is not available in the common theory.  Further, we have seen from the discussion in Section~\ref{sec:pwnot} that the definition of $(\M_k)^{-1}$ has only been implied from definitions of $\M_k$.  In this section, we define a set of interpolation functions $\ol\I$ analogous to the Whitney functions and use them to provide an explicit definition of a dual discrete Hodge star.  

\begin{figure}[ht]
\psfrag{vi}{$\bv_i$}
\psfrag{Ci}{$C_i$}
\psfrag{x}{$\bx$}
\psfrag{Dix}{$D(\bx)$}
\psfrag{DixCi}{$D(\bx)\cap C_i$}
\[\begin{array}{cc}
\includegraphics[width=.3\linewidth]{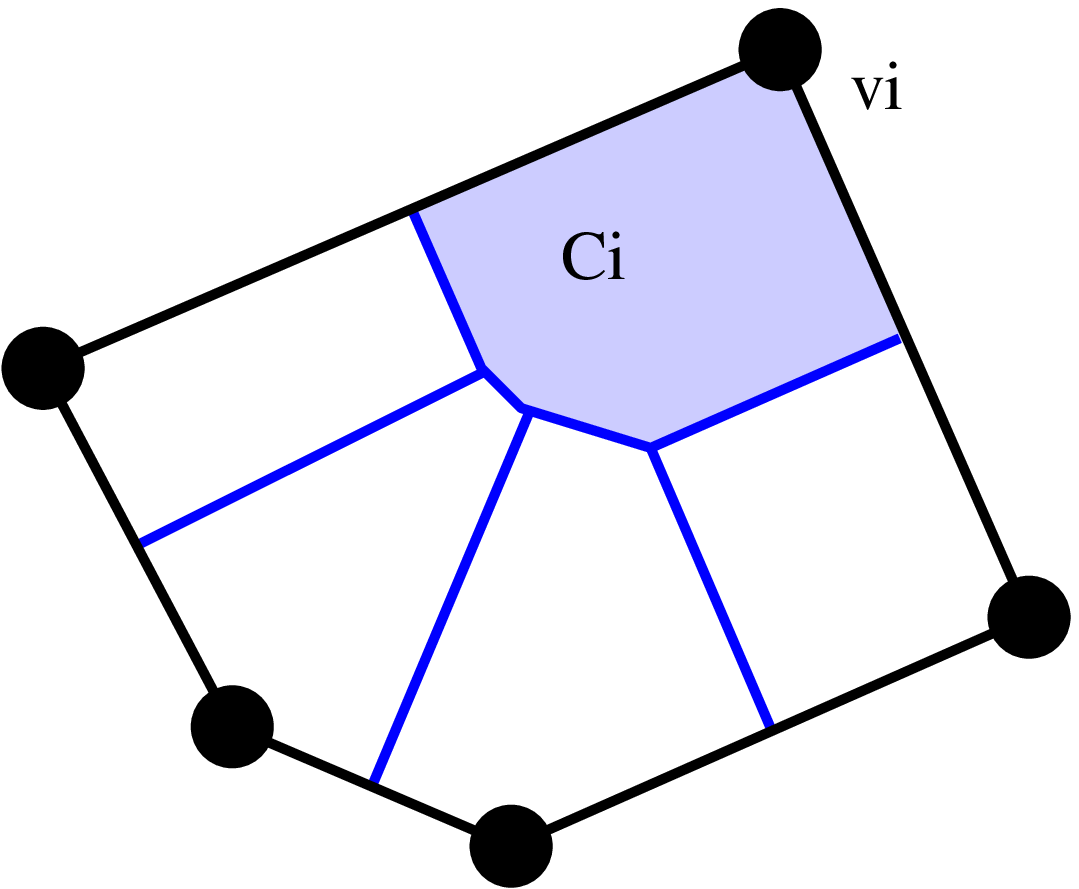} &
\includegraphics[width=.3\linewidth]{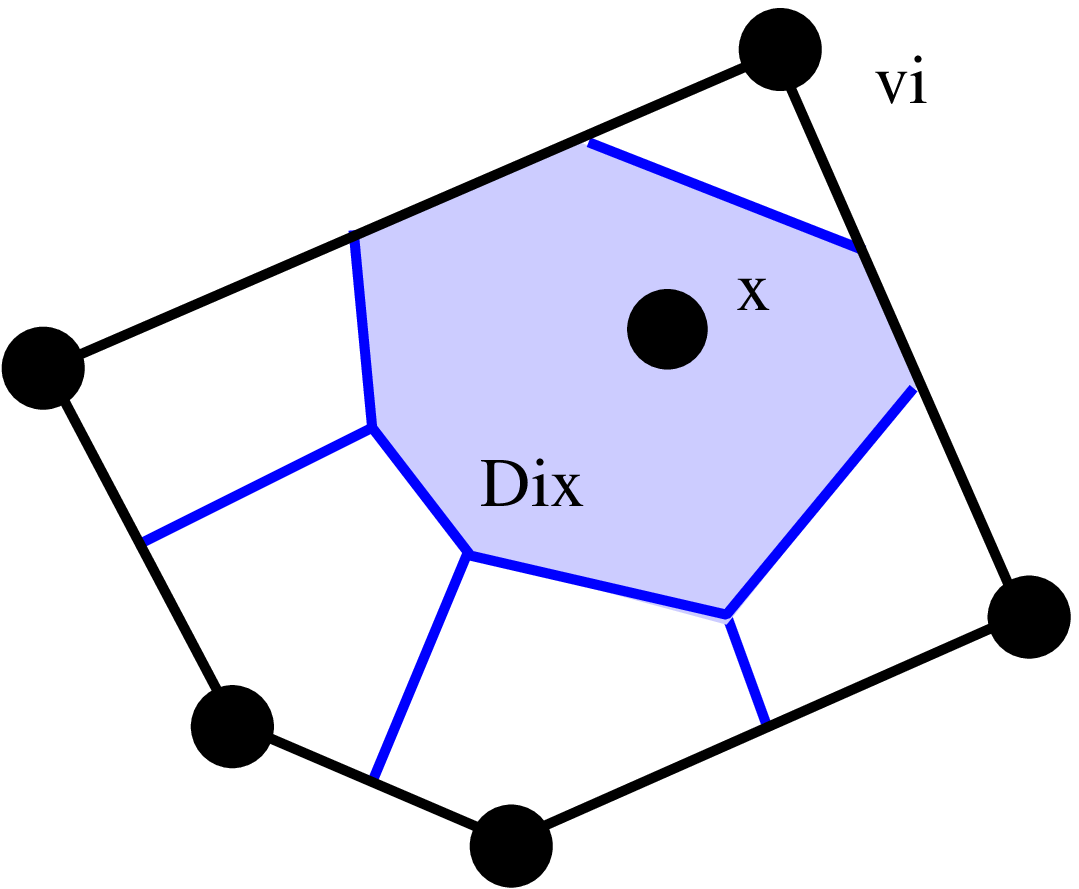} \\
\includegraphics[width=.3\linewidth]{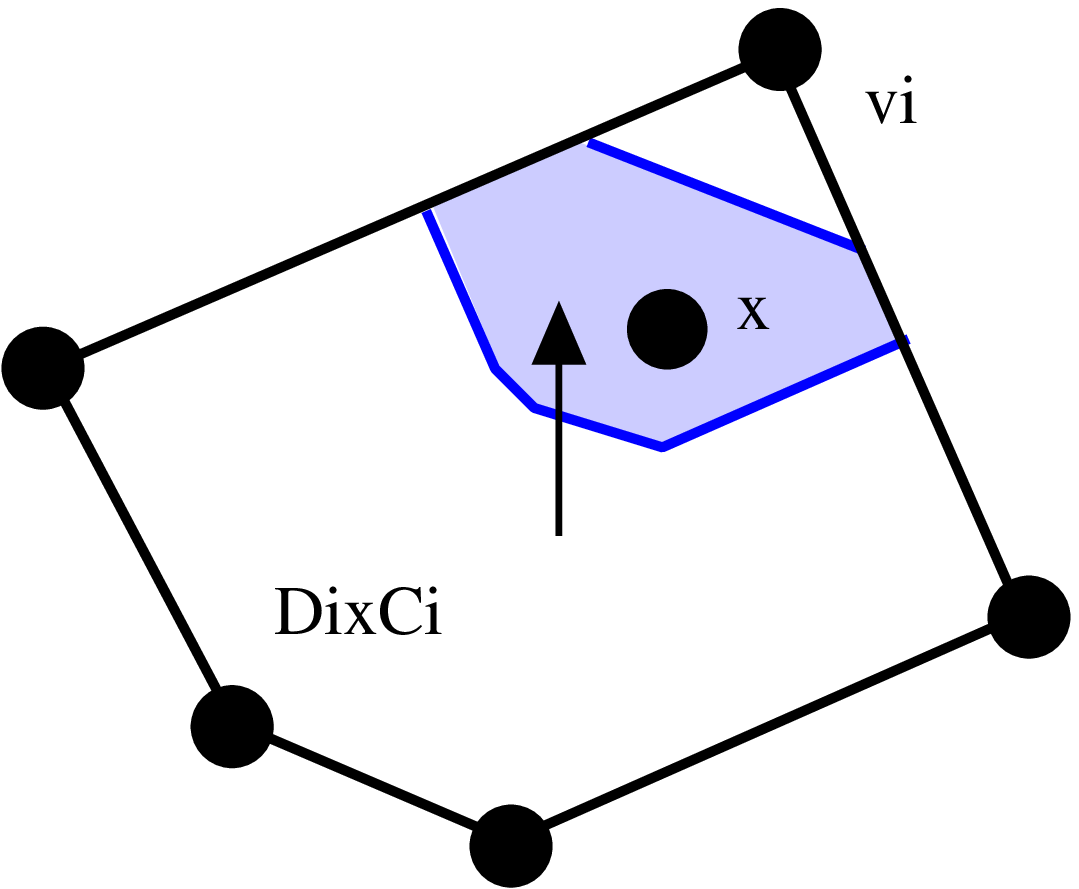} &
\includegraphics[width=.3\linewidth]{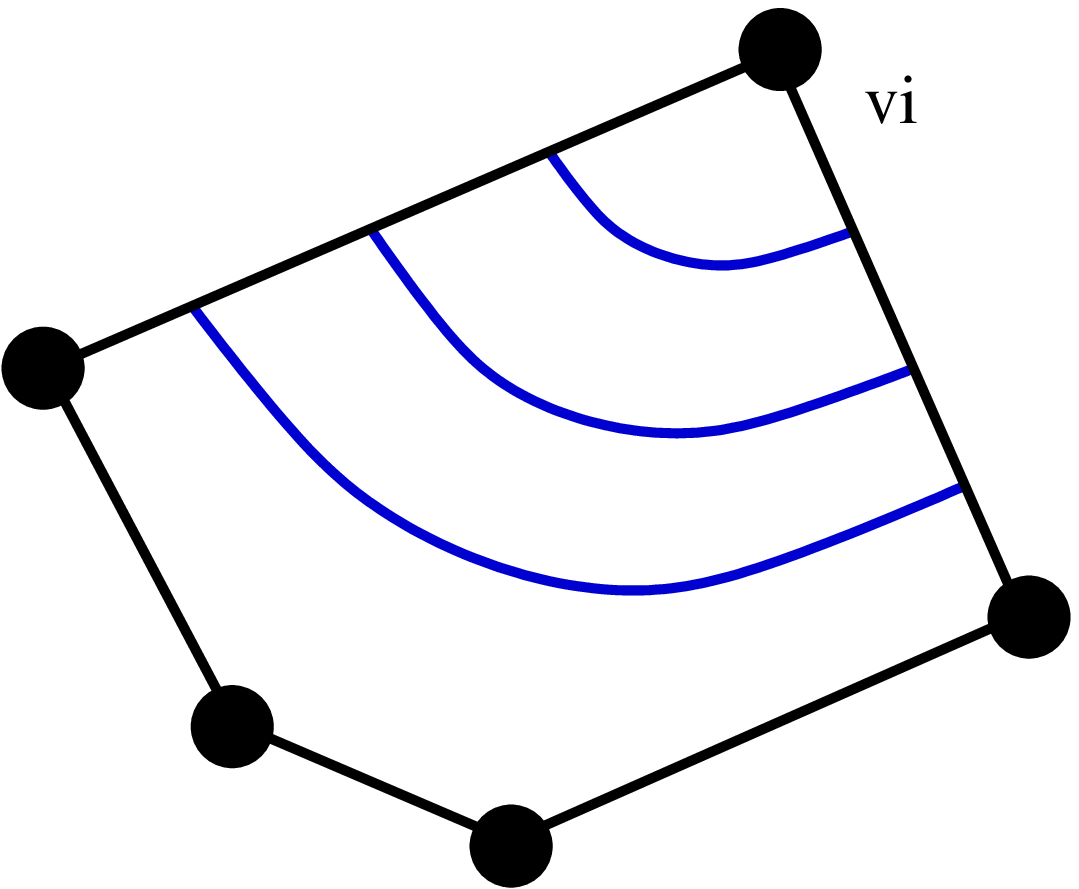}
\end{array}\]
\caption{Geometric calculation of a Sibson coordinate.  $C_i$ is the area of the Voronoi region associated to vertex $\bv_i$ inside $\Tau$.  $D(\bx)$ is the area of the Voronoi region associated to $\bx$ if it is added to the vertex list.  The quantity $D(\bx)\cap C_i$ is exactly $D(\bx)$ if $\bx=\bv_i$ and decays to zero as $\bx$ moves away from $\bv_i$, with value identically zero at all vertices besides $\bv_i$.  The bottom right figure shows how the level sets of the Sibson coordinate associated to $\bv_i$ sit inside a single polygon.  More figures can be found in Milbradt and Pick~\cite{MP2008}}
\label{fig:sibson}
\end{figure}

Define the dual Whitney interpolant of a dual $k$-cochain $\dc w\in\dcs^{k}$ to be
\begin{equation}
\label{eqn:olikdef}
\ol\I_{k}(\dc w):=\sum_{{\star\sigma^{n-k}}\in\dcs_{k}}\dc w({\star\sigma^{n-k}})\dwhit_{\star\sigma^{n-k}}
\end{equation}
where $\dwhit_{\star\sigma^{n-k}}$ is a dual Whitney function associated to the $k$-cell $\star\sigma^{n-k}$ in the dual mesh.  These functions are defined using a generalization of barycentric coordinates known as Sibson functions~\cite{S1980}, also called the natural neighbor or natural element coordinates~\cite{SM2006}.  Figure~\ref{fig:sibson} summarizes the definition.

\begin{defn}
\label{def:sibsonfns}
{\em
Let $\bx$ be a point inside a polyhedral cell $\Tau$ of the dual mesh.  Let $P$ denote the set of vertices $\{\bv_i\}$ and define
\[P'=P\cup \{\bx\} = \{\bv_1,\ldots,\bv_{N},\bx\}.\]
Denote the \textbf{Voronoi cell} associated to a point $\bp$ in a pointset $Q$ by
\[V_Q(\bp):= \left\{\by\in \Tau \, : \, \vsn{\by-\bp} < \vsn{\by-\bq} \, , \, \forall\bq\in Q\setminus\{\bp\}  \right\}.\]
Note that these Voronoi cells have been restricted to $\Tau$ and are thus always of finite size.  Fix the notation
\[\begin{array}{lcl}
C_i &:=& |V_{P}(\bv_i)| = \left|\{\by\in\Tau \, : \, \vsn{\by-\bv_i} < \vsn{\by-\bv_j}\, , \, \forall j\not=i\}\right| \\
 & =& \text{area of cell for $\bv_i$ in Voronoi diagram on the points of $P$,} \\
\\
D(\bx) &:=&  |V_{P'}(\bx)| = \left|\{\by\in\Tau \, : \, \vsn{\by-\bx} < \vsn{\by-\bv_i}\, , \, \forall i\}\right| \\
 & =& \text{area of cell for $\bx$ in Voronoi diagram on the points of $P'$}.
\end{array}
\]
By a slight abuse of notation, define
\[ D(\bx)\cap C_i := |V_{P'}(\bx)\cap V_{P}(\bv_i)|. \]
The notation is shown in Figure~\ref{fig:sibson}.  The Sibson coordinates are defined to be
\begin{align*}
\ol\lambda_i(\bx) &:= \frac{D(\bx)\cap C_i}{D(\bx)} &
\textnormal{or, equivalently,} & &
\ol\lambda_i(\bx) &= \frac{D(\bx)\cap C_i}{\sum_{j=1}^{N}D_j(\bx)\cap C_j}.
\end{align*}
}
\end{defn}

Milbradt and Pick \cite{MP2008} modify the definition of the Sibson functions for polytopes so that the coordinates of a point on an edge or facet of the polytope are dependent only on the Sibson functions associated to the boundary vertices of that edge or facet.  This ensures $C^0$ continuity of the functions across adjacent mesh elements.

Moreover, it has been shown that the Sibson functions are $C^\infty$ on the polygon except at the vertices $\bv_i$ where they are $C^0$ and on circumcircles of Delaunay triangles where they are $C^1$~\cite{S1980,F1990}.  Since the finite set of vertices are the only points at which the function is not $C^1$, we conclude that $\ol\lambda_i\in H^1(K)$ where $K$ is the domain mesh.  This is the typical continuity required for finite element applications with nodal interpolation functions and makes them fit for use in the dual Whitney functions we define next.

\begin{figure}[ht]
\centering
\[\begin{array}{ccc}
\includegraphics[height=2.5in]{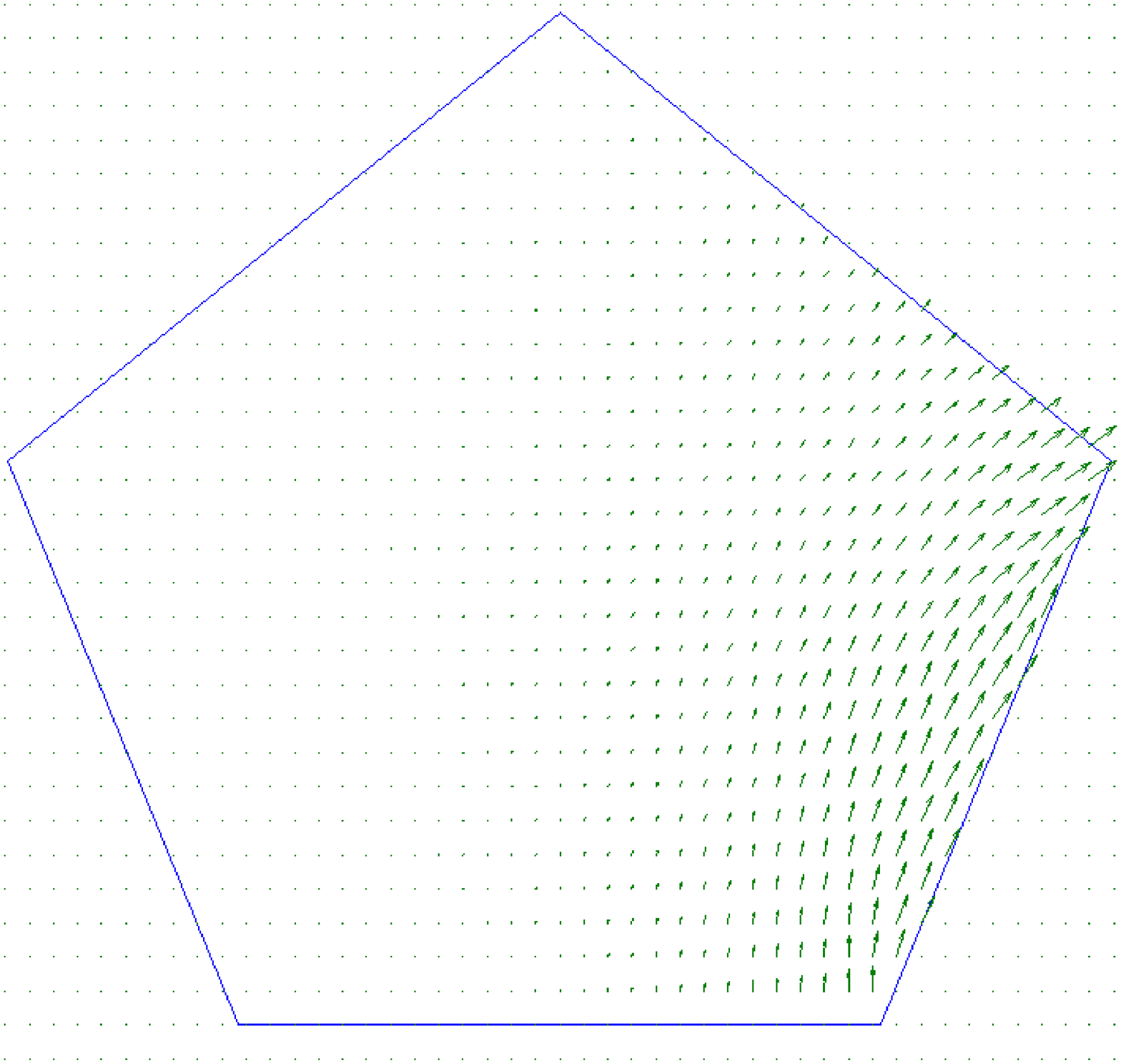} & \quad &
\quad\includegraphics[height=2.5in]{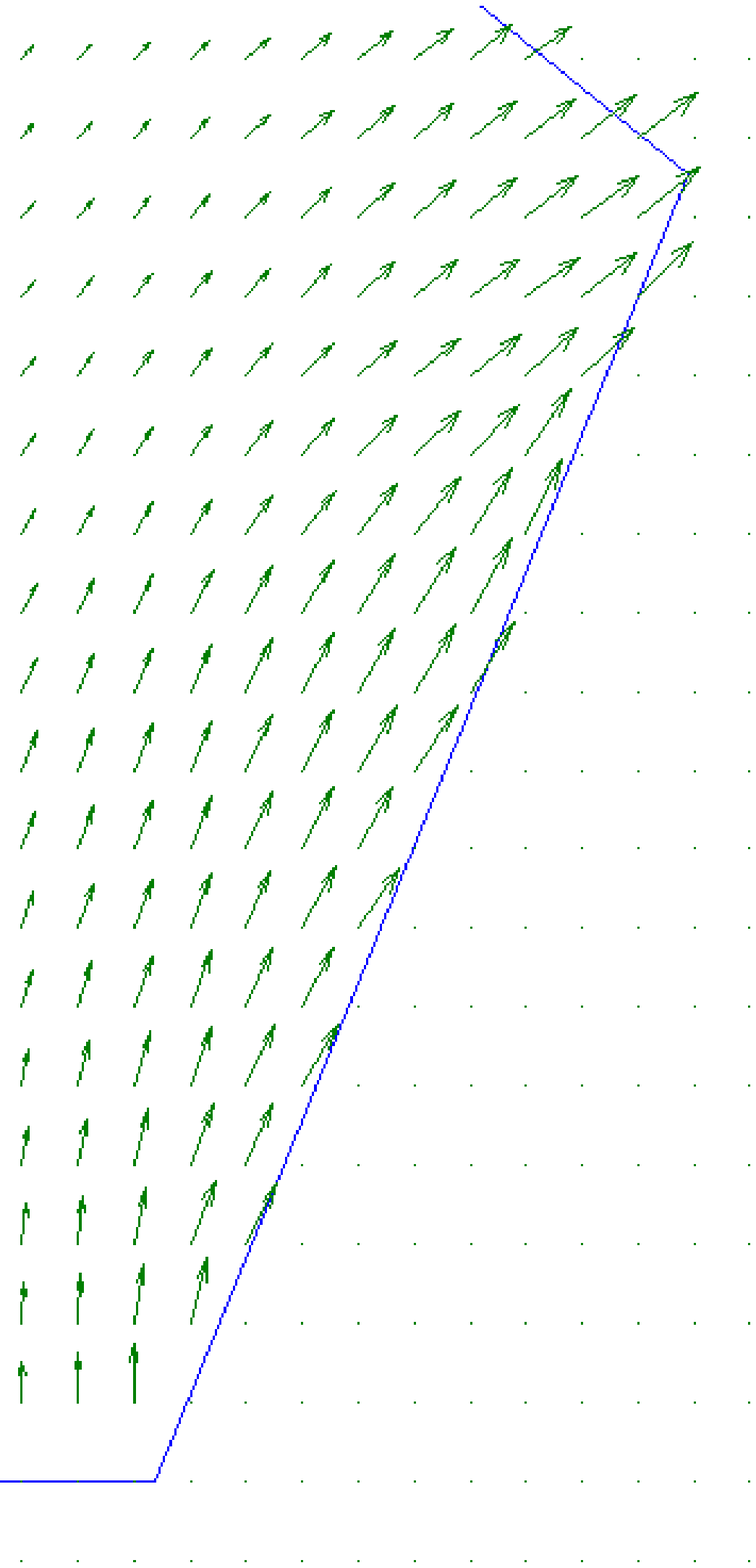}
\end{array}
\]
\caption{The dual Whitney function associated to the lower right edge of a pentagon is shown on the left.  The magnified portion shows the vector field in the neighborhood of this edge.  The gradients were approximated in Matlab using a simple 2-point difference rule on a regular grid laid over the pentagon. }
\label{fig:whitDualOneForm}
\end{figure}

\begin{figure}[ht]
\centering
\psfrag{v1}{$\ol\bv_1$}
\psfrag{v2}{$\ol\bv_2$}
\psfrag{v3}{$\ol\bv_3$}
\psfrag{v4}{$\ol\bv_4$}
\psfrag{v0}{$\ol\bv_0$}
\psfrag{c}{$\ol\bc$}
\psfrag{T0}{$\tau_0$}
\psfrag{T1}{$\tau_1$}
\psfrag{T2}{$\tau_2$}
\psfrag{T3}{$\tau_3$}
\psfrag{T4}{$\tau_4$}
\psfrag{s1}{\textcolor{DkBlue}{$\sigma^1$}}
\[\begin{array}{cc}
\includegraphics[width=.35\linewidth]{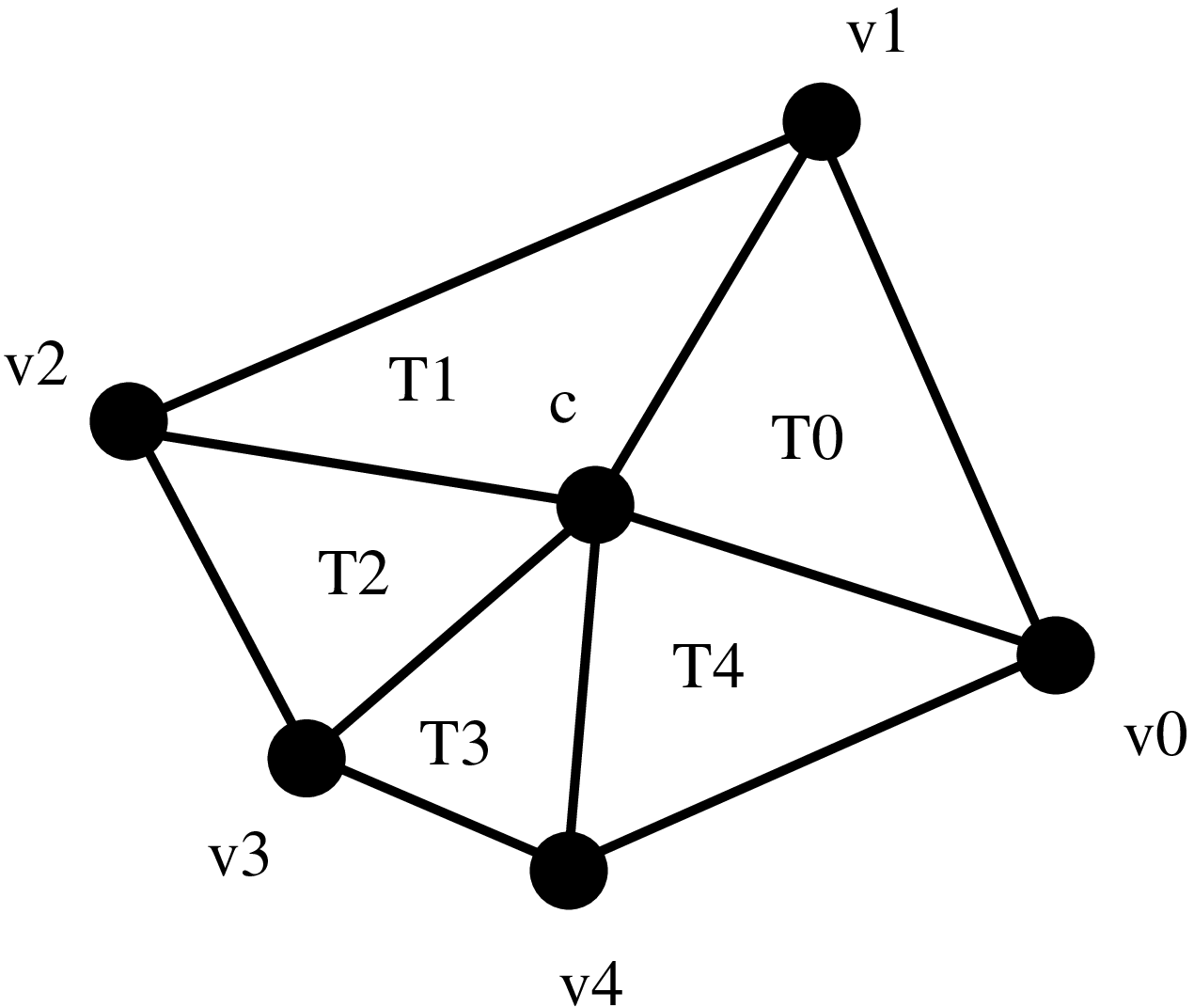} &
\quad\includegraphics[width=.35\linewidth]{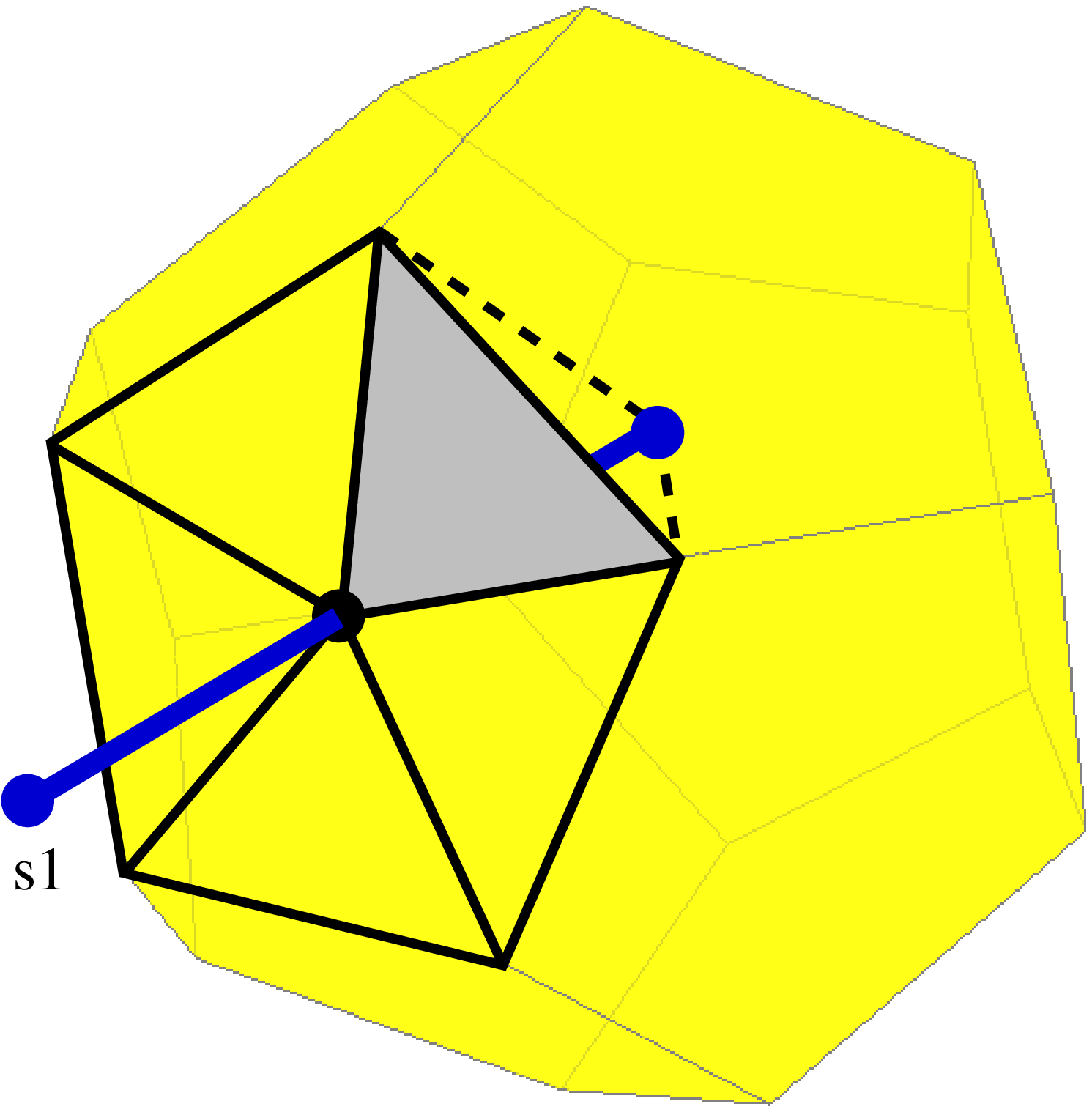}
\end{array}
\]
\caption{Sample computation of a dual Whitney function associated to a dual face $\star\sigma^1$ with vertices $\ol\bv_i$. By adding the centroid $\ol\bc$, we have a canonical decomposition of $\star\sigma^1$ into triangles $\tau_i$.  A weighted sum of the primal Whitney function associated with each $\tau_i$ is constructed to define the function for the face.  As shown on the right, each $\tau_i$, e.g. the shaded triangle, forms a tetrahedron by connecting its vertices to the vertex of $\sigma^1$ interior to the polyhedron.  Note that in general $\ol\bc$ need not be the same as $\sigma^1\cap\star\sigma^1$. }

\label{fig:whitDualTwoForm}
\end{figure}

\begin{defn}
\label{def:whitlikefns3d}
{\em
The \textbf{dual Whitney function} $\dwhit_{\star\sigma^{3-k}}$ associated to the $k$-dimensional element $\star\sigma^{3-k}$ in a 3D dual mesh is defined as follows.
\begin{itemize}
\item\textbf{Dual Vertices}.  The function associated to a dual vertex $\star\sigma^3 :=\ol\bv_i$ is the Sibson coordinate for the vertex, i.e.
\[\dwhit_{\star\sigma^3}:=\ol\lambda_i\]
\item\textbf{Dual Edges}.  The function associated to an oriented dual edge $\star\sigma^2:=[\ol\bv_i,\ol\bv_j]$ is the vector-valued function
\[\dwhit_{\star\sigma^2}:=\ol\lambda_i\nabla\ol\lambda_j-\ol\lambda_j\nabla\ol\lambda_i\]
An example is shown in Figure~\ref{fig:whitDualOneForm}.

\item\textbf{Dual Faces}.  Consider a dual face $\star\sigma^1$ with $m$ vertices $\{\ol\bv_0,\ldots,\ol\bv_{m-1}\}$.  Partition the face canonically into triangles by adding a vertex $\ol\bc$ at the centroid of the face vertices and adding the edges $[\ol\bc,\ol\bv_i]$.  Define 2-simplices $\tau_i := [\ol\bc,\ol\bv_i,\ol\bv_{i+1}]$, indices taken mod $m$.  Define 3-simplices by connecting the $\tau_i$ to the endpoint of $\sigma^1$ inside the polyhedron.  Define
\[\dwhit_{\star\sigma^1} := \ds \sum_{i=0}^{m-1}\frac{|\tau_i|}{|\star\sigma^1|}\whit_{\tau_i}\chi_{\tau_i},\]
where $\chi_{\tau_i}$ is the characteristic function on $\tau_i$ (1 on $\tau_i$, 0 otherwise) and
\[\whit_{\tau_i} := 2\left(\whitC {\ol\bc} i {i+1} -\whitC i {\ol\bc} {i+1} + \whitC {i+1} {\ol\bc} i \right).\]
Note that $\whit_{\tau_i}$ is the Whitney 2-form associated to face $\tau_i$ of a tetrahedron (see (\ref{eq:prwhit2}) in Appendix A) and that these tetrahedra partition the entire polyhedra.  An example is shown in Figure~\ref{fig:whitDualTwoForm}.
\item\textbf{Dual Cells}.   The scalar-valued function associated to a dual cell $\star\sigma^0$ is a constant function on the cell:
\[\dwhit_{\star\sigma^0}:= \chi_{\star\sigma^0} = \left\{\begin{array}{rl} 1/|\star\sigma^0| & \text{on $\star\sigma^0$}\\ 0 & \text{otherwise}  \end{array}\right.  \]
\end{itemize}
}
\end{defn}

Since the dual Whitney functions use a generalization of barycentric coordinates, it can be shown that they have the standard continuity across faces, e.g. tangential continuity for $\dwhit_{\star\sigma^2}$  and normal continuity for $\dwhit_{\star\sigma^1}$.  This means the image of $\ol\I_0$ is in $H^1$, the image of $\ol\I_1$ is in $\hcurl$, and so forth (see Figure~\ref{fig:derham}).  A proof of this and other properties of $\ol\I_k$ appears in~\cite{G2011}.  We are also developing a higher order version of these operators~\cite{GB2011m}.

Using dual Whitney functions, we define a novel dual discrete Hodge star by
\begin{equation}
\label{eqn:dhswhit}
((\M_k^{Dual})^{-1})_{ij} := \left(\dwhit_{\star\sigma^k_i},\dwhit_{\star\sigma^k_j}\right).
\end{equation}
The inner product here is the standard integration of scalar or vector valued functions over the dual domain $\star K$. For instance, in the case $k=3$, the definition yields
\[((\M_3^{Dual})^{-1})_{ij} := \left(\dwhit_{\star\sigma^3_i},\dwhit_{\star\sigma^3_j}\right) = \int_{\star K} \ol\lambda_i\ol\lambda_j.\]
The formulation for other $k$ values will similarly involve integrals of the $\ol\lambda_i$ functions.
\begin{lemma}
\label{lem:mkdsparse}
$(\M_k^{Dual})^{-1}$ is sparse.
\end{lemma}
\begin{proof}
Observe that $\dwhit_{\star\sigma^k}$ has localized support by construction.  Entry $ij$ of $(\M_k^{Dual})^{-1}$ will be non-zero only if $\star\sigma_i^k$ and $\star\sigma_j^k$ are adjacent.  Thus each row of the matrix will have at most as many non-zero entries as $\star\sigma^k_i$ has adjacent $n-k$ cells, meaning the matrix is sparse.
\end{proof}
Lemma \ref{lem:mkdsparse} does not hold if $\M_k^{Dual}$ is replaced by $\M_k^{Whit}$ as these sparse matrices typically have dense inverses.  Note that $(\M_k^{Diag})^{-1}$ is trivially sparse since it is diagonal, however, it can only be employed when the meshes are orthogonal.

\subsection{Local Structure of Discrete Hodge Stars}
The continuous Hodge star $\ast$ is a local operator meaning its effect on a differential form evaluated at a particular point on a manifold depends only on the geometry of a local neighborhood of the point.  The discrete Hodge star is thus required to be a local operator as well meaning the evaluation of $\M_k$ on a basis cochain $\pc w^k_i$ (1 on $\sigma_i^k$ and 0 otherwise) should involve values on only a few simplices adjacent to $\sigma_i^k$.  In the language of matrix theory, this requirement says $\M_k$ should be sparse.

\begin{figure}[ht]
\[\begin{array}{ccc}
\includegraphics[height=1.5in]{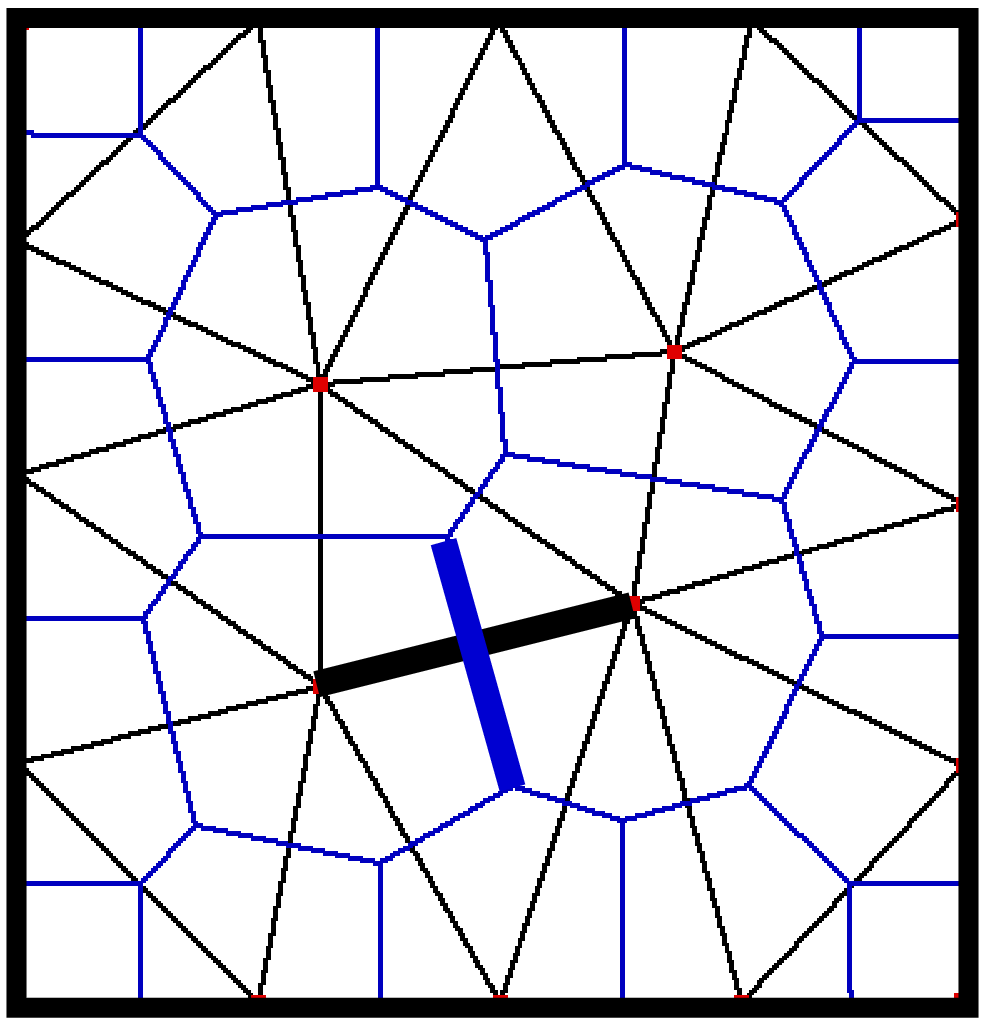} &
\includegraphics[height=1.5in]{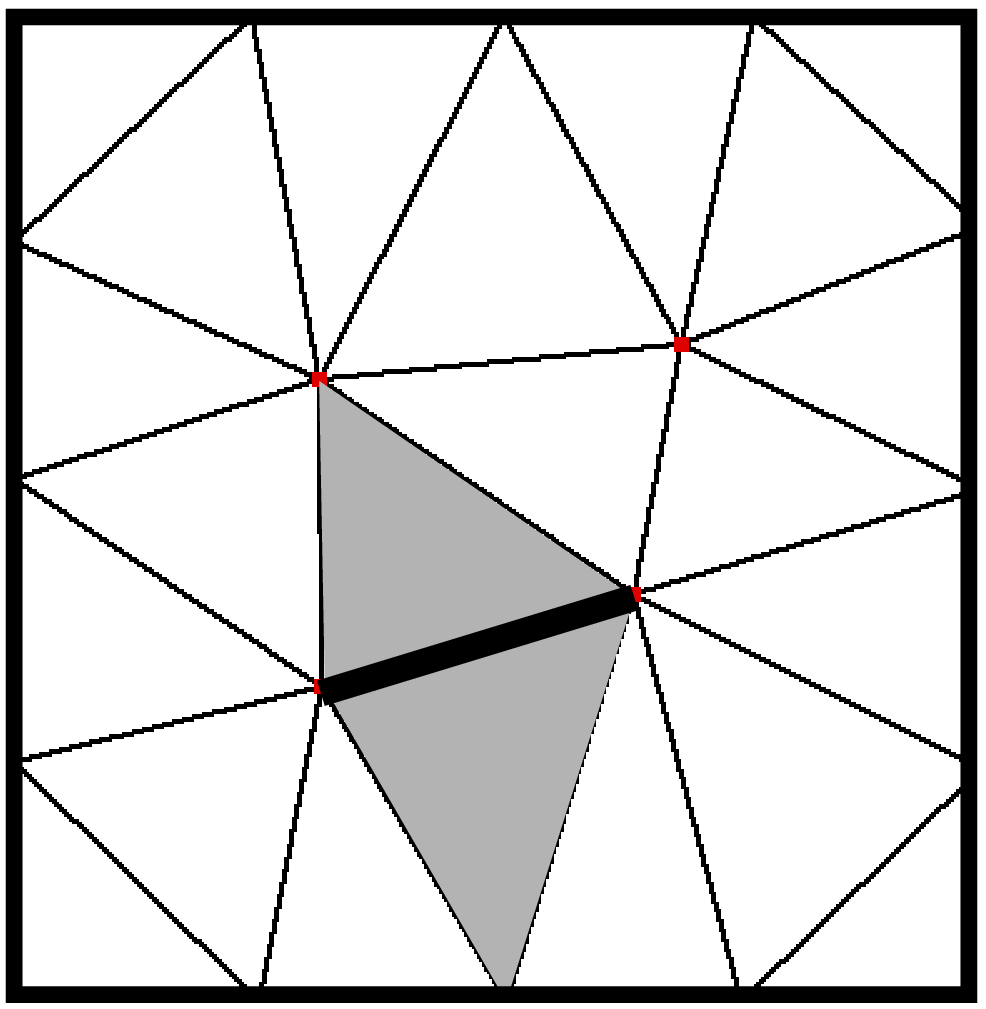} &
\includegraphics[height=1.5in]{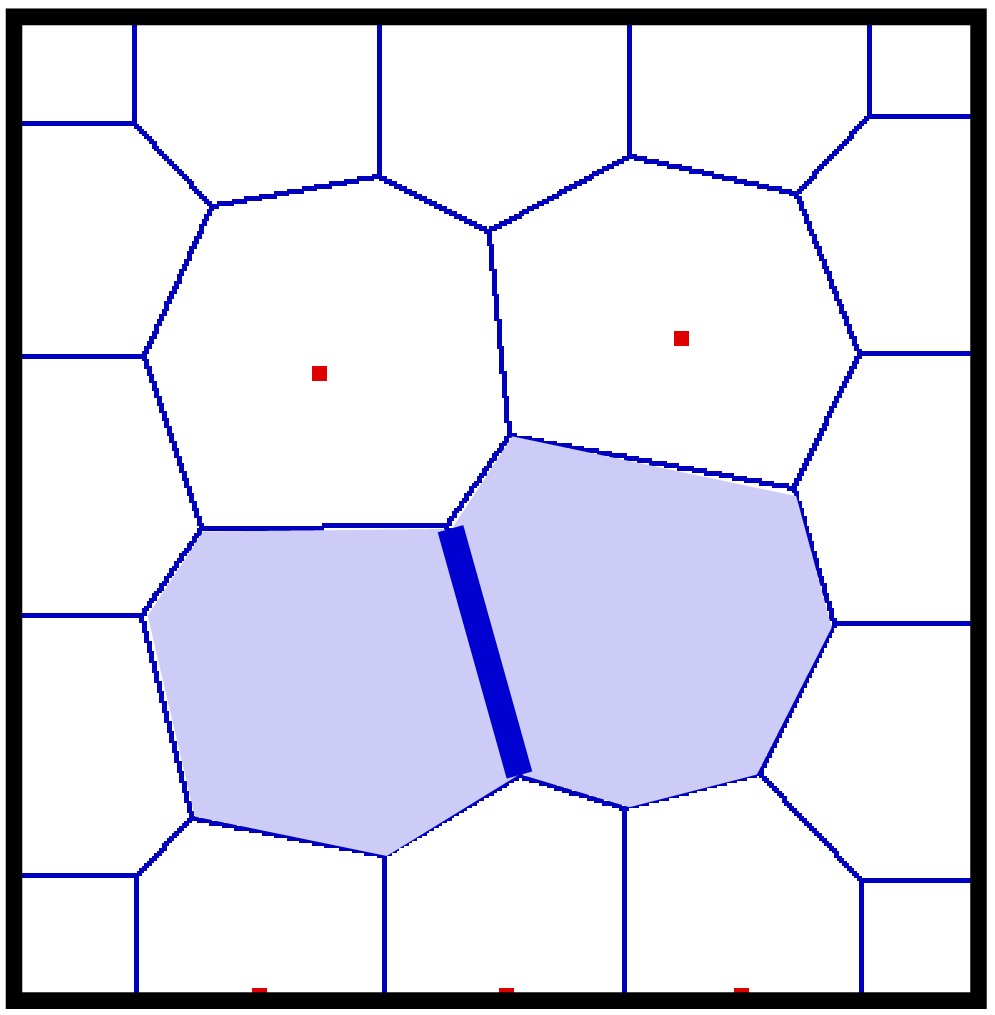} \\
\Mdiag 1 & \Mwhit 1 & \Mduali 1
\end{array}\]
\caption{The various discrete Hodge stars depend on different aspects of mesh geometry as shown in this 2D examples.  The diagonal Hodge star (left) computes ratios of sizes of primal-dual element pairs.  The Whitney Hodge star (middle) has entries of Whitney functions integrated against each other.  The support of a particular $\whit_{\sigma^1_i}$ function is shown in grey; the integral of its projection to the bold edge has value 1.  The Dual Hodge star (right) that we propose has entries of dual Whitney functions integrated against each other.  The support of a particular $\dwhit_{\star\sigma^1_i}$ is shown in blue; the integral of its projection to the bold dual edge has value 1.}
\label{fig:dhstypes}
\end{figure}

We now give a more specific characterization of the sparsity structure of $\M_k^{Whit}$ and $(\M_k^{Dual})^{-1}$.  The intuition for these results is demonstrated by Figure~\ref{fig:dhstypes}

\begin{lemma}
Entry $ij$ in $\M_k^{Whit}$ is non-zero only if there exists $\sigma^n\in K$ such that $\sigma^n$ has at least one vertex from $\sigma^k_i$ and one vertex from $\sigma^k_j$.
\end{lemma}
\begin{proof}
Computing entry $ij$ in $\M_k^{Whit}$ involves \cite[Prop. 9.6]{BH2011} summing terms of the form
\begin{equation}
\label{eqn:hswhitcomp}
\left(\int_K\lambda_1\lambda_2\right)\det\left(V_I^T W_J\right)
\end{equation}
where $\lambda_1,\lambda_2$ are barycentric functions associated to $v_1\in\sigma^k_i$, $v_2\in\sigma^k_j$, respectively; $I$ is a list of $k$ vertices from $\sigma^k_i$ not including $v_1$; $J$ is a list of $k$ vertices from $\sigma^k_j$ not including $v_2$; and $V_I$, $W_J$ are $n\times k$ matrices.  The $p$th column of $V_I$ is the vector $\nabla\lambda_p$ where $\lambda_p$ is the barycentric function associated to the $p$th entry in $I$.  The $q$th column of $W_J$ is the vector $\nabla\lambda_q$ where $\lambda_q$ is the barycentric function associated to the $q$th entry in $J$.

Observe that the support of the barycentric function associated to vertex $v$ is contained within the $n$-simplices touching $v$.  Thus, if there is no $\sigma^n$ with at least one vertex from $\sigma^k_i$ and one vertex from $\sigma^k_j$, the $\lambda_1$ and $\lambda_2$ appearing in (\ref{eqn:hswhitcomp}) will always have disjoint support, making the entry zero.
\end{proof}
Using the same kind of reasoning, we have a similar result for our dual discrete Hodge star.
\begin{lemma}
Entry $ij$ in $(\M_k^{Dual})^{-1}$ is non-zero only if there exists $\star\sigma^0\in \star K$ such that $\star\sigma^0$ has at least one vertex from $\star\sigma^k_i$ and one vertex from $\star\sigma^k_j$.
\end{lemma}
The number of $k$-simplices in an $n$-simplex is ${n+1\choose k+1}$ which gives the following corollary.

\begin{cor}
Let $A(\sigma^k)$ denote the number of $n$-simplices in $K$ incident on at least one vertex from $\sigma^k$.  Then the number of non-zero entries in row $i$ of $\M_k^{Whit}$ or row $i$ of $(\M_k^{Dual})^{-1}$ is at most ${n+1\choose k+1} A(\sigma^k_i)$.
\end{cor}
The bound can be sharpened for particular choices of $n$ and $k$ or if additional assumptions are made about $K$.  As stated, however, the corollary provides a simple means for evaluating the computational expense of a particular discretization scheme as we will discuss in Section \ref{sec:res}.

\subsection{Numerical Stability}
To maintain the numerical stability of a DEC-based method, the discrete Hodge star matrix should have a bounded condition number.  Put differently, the entries of the matrix should be roughly the same order of magnitude.  This requirement is frequently considered from the context of numerical analysis but is often absent from the literature on discrete operators.

\begin{figure}[ht]
\centering
\psfrag{(a)}{\textbf{(a)}}
\psfrag{(b)}{\textbf{(b)}}
\psfrag{(c)}{\textbf{(c)}}
\psfrag{sigma1}{$\sigma^1$}
\psfrag{starsigma1}{\textcolor{red}{$\star\sigma^1$}}
\psfrag{sigma2}{$\sigma^2$}
\psfrag{starsigma2}{\textcolor{red}{$\star\sigma^2$}}
\[\begin{array}{ccc}
\includegraphics[width=.3\linewidth]{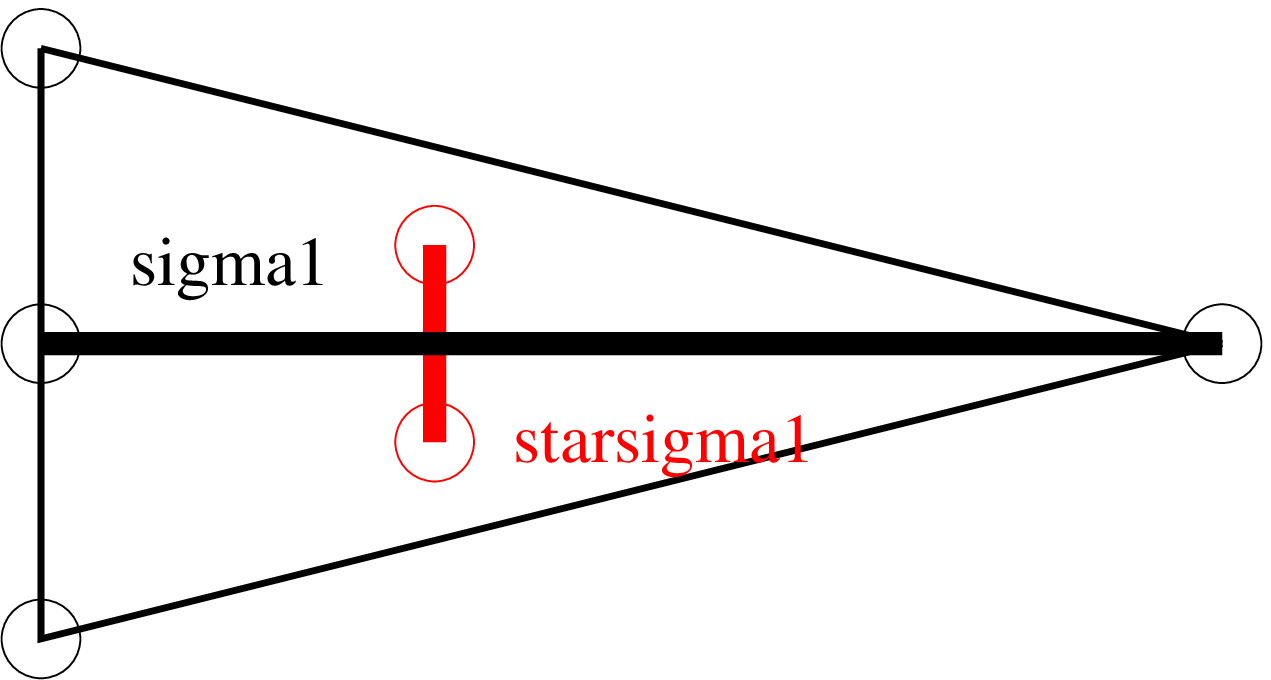} &
\includegraphics[width=.3\linewidth]{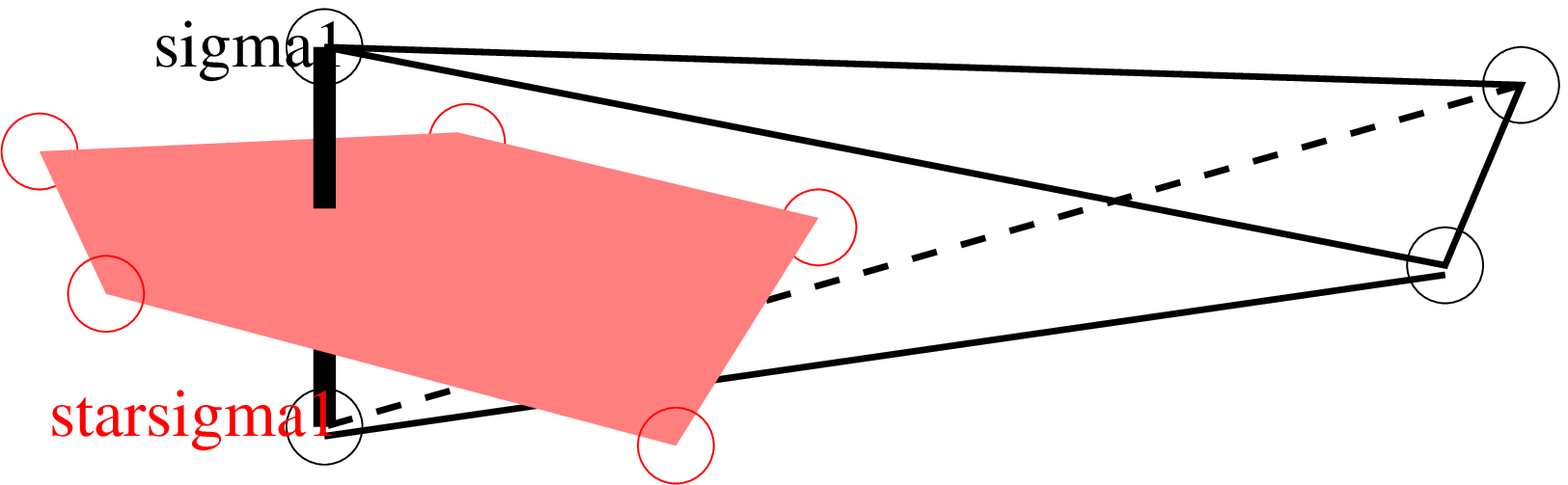} &
\includegraphics[width=.3\linewidth]{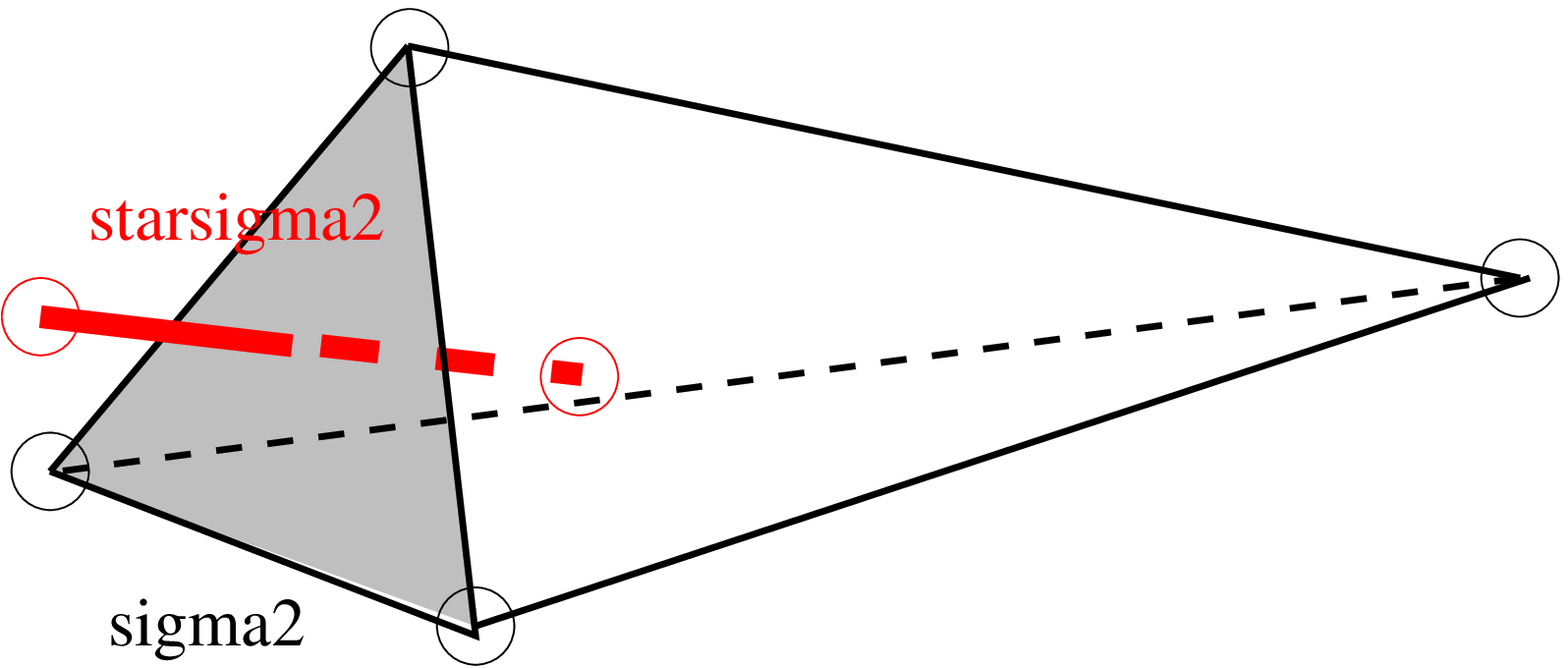} \\
\textbf{(a)} & \textbf{(b)}  &  \textbf{(c)}
\end{array}\]
\caption{Examples illustrating how the measure of a primal simplex $\sigma^k$ (black) and its dual $\star\sigma^k$ (red) need not be the same order of magnitude.  \textbf{(a)} In this 2D example, the ratio $|\star\sigma^1|/|\sigma^1|$ can be made arbitrarily small by increasing the length of $\sigma^1$.  \textbf{(b)}  The ratio $|\star\sigma^1|/|\sigma^1|$ can be made arbitrarily large by decreasing the length of $\sigma^1$.  \textbf{(c)} The ratio $|\star\sigma^2|/|\sigma^2|$ can be made arbitrarily large by decreasing the area of $\sigma^2$.  Thus, a discrete Hodge star involving terms of the form $|\star\sigma^k|/|\sigma^k|$ may have a bad condition number unless primal \textit{and} dual mesh quality is controlled.}
\label{fig:dualratioprobs}
\end{figure}

\begin{figure}[ht]
\centering
\includegraphics[width=.3\linewidth]{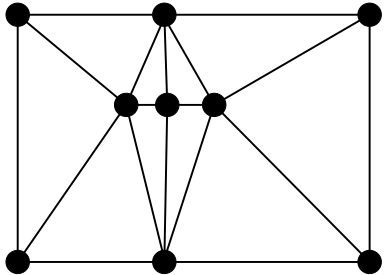}
\caption{Graded meshes also present a problem for discrete Hodge stars involving primal-dual size ratios. The primal mesh shown here induces a wide variation in values of \newline $|\star\sigma^k|/|\sigma^k|$ for $k=0,1,2$.  This can cause ill-conditioned $\M_k$ matrices, resulting in numerical instability.}
\label{fig:gradationprob}
\end{figure}

The common thread in the geometrically-defined discrete Hodge stars such as $\M_k^{Diag}$ is a measurement of the size of dual cells i.e. $|\star\sigma^k|$.  This suggests that geometric criteria on primal elements alone will not be sufficient to control the condition number of the discrete Hodge star matrix.  In particular, since ratios of primal to dual cells are computed, the following criteria must be satisfied:
\renewcommand{\labelenumi}{N\arabic{enumi}.}
\begin{enumerate}
\setcounter{enumi}{0}
\item Primal simplices $\sigma^k$ satisfy geometric quality measures.
\item Dual cells $\star\sigma^k$ satisfy geometric quality measures.
\item The value of $|\star\sigma^k|/|\sigma^k|$ is bounded above and below.
\item The primal and dual meshes do not have large gradation of elements, i.e. $\min_i |\sigma_i^k|$ and $\max_i |\sigma_i^k|$ are the same order of magnitude and $\min_i |\star\sigma_i^k|$ and $\max_i |\sigma_i^k|$ are the same order of magnitude.
\end{enumerate}
Conditions N1 and N2 are required for discretization stability.  Aspect ratio is often used as a geometric quality measure for tetrahedra.  Conditions N3 and N4 are based on our analysis above.  Condition N4 in particular shows that these discrete Hodge stars are not fit for use on meshes tailored to multi-resolution situations where gradation is necessary to achieve reasonable computation times.  Examples are shown in Figures \ref{fig:dualratioprobs} and \ref{fig:gradationprob}.

For $\M_k^{Whit}$, the size of the matrix entries are controlled by the size of the inner products of Whitney basis forms.  The integrals in (\ref{eqn:hswhitcomp}) are on the order of the size of $|\sigma_k|$, meaning again that a large gradation in primal mesh element size could produce large condition numbers.  Since $\M_k^{Whit}$ does not depend on the size of dual mesh elements, however, its condition number is more stable against violations of conditions N2 and N3.  Analogously, the condition number of $(\M_k^{Dual})^{-1}$ is more stable against violations of conditions N1 and N3.  Our conclusions are summarized below.

\begin{itemize}
\item Conditions N1-N4 are necessary to ensure $\M_k^{Diag}$ has a good condition number.
\item Conditions N1 and N4 are necessary to ensure $\M_k^{Whit}$ has a good condition number.
\item Conditions N2 and N4 are necessary to ensure $(\M_k^{Dual})^{-1}$ has a good condition number.
\end{itemize}

\subsection{Improved Condition Numbers with $(\M_k^{Dual})^{-1}$}

To provide concrete evidence for our numerical stability claims, we present a simple example in 2D showing how $\M_1^{Diag}$ and $\M_1^{Whit}$ can have condition numbers an order of magnitude worse than $(\M_1^{Dual})^{-1}$ on the same mesh.  This serves as a proof of concept that the DEC-based dual formulation of a problem can provide practical advantages in cases of difficult mesh geometry.

\begin{figure}[ht]
\centering
\psfrag{sigma12}{$\sigma_{12}$}
\psfrag{sigma13}{$\sigma_{13}$}
\psfrag{sigma14}{$\sigma_{14}$}
\psfrag{sigma23}{$\sigma_{23}$}
\psfrag{sigma24}{$\sigma_{24}$}
\psfrag{v1}{$\bv_1$}
\psfrag{v2}{$\bv_2$}
\psfrag{v3}{$\bv_3$}
\psfrag{v4}{$\bv_4$}
\includegraphics[width=.4\linewidth]{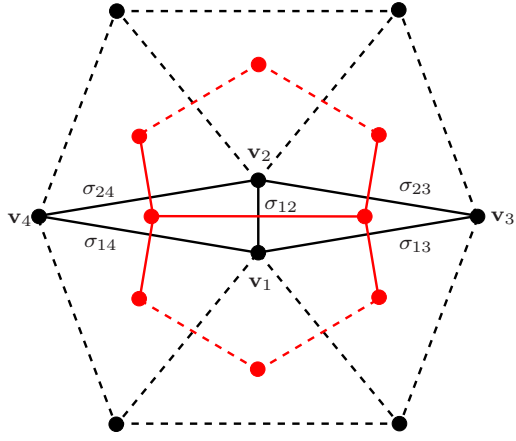}
\caption{Mesh used for sample calculation of $\M_1$ matrices.  The vertices have coordinates $\bv_1=(0,0)$, $\bv_2=(0,1)$, $\bv_3=(P,\frac 12)$, $\bv_4=(-P,\frac 12)$.}
\label{fig:M1compNotn}
\end{figure}

In the 2D mesh shown in Figure~\ref{fig:M1compNotn}, the labeled vertices of the primal mesh have coordinates $\bv_1=(0,0)$, $\bv_2=(0,1)$, $\bv_3=(P,\frac 12)$, and $\bv_4=(-P,\frac 12)$, where $P$ is a free parameter we can adjust to modify the geometry.  The remaining vertices are chosen so that they form equilateral triangles with edges $\sigma_{13}$, $\sigma_{23}$, $\sigma_{14}$, and $\sigma_{24}$, as shown.  The orthogonal, circumcenter-based dual mesh is shown in red.

Without loss of generality, fix any ordering on the mesh edges, beginning with
\begin{equation}
\label{eq:edgeorder}
\{\sigma_{12},\sigma_{13},\sigma_{14},\sigma_{23},\sigma_{24},\ldots\}.
\end{equation}
We first calculate the upper left $5\times 5$ block of each matrix, yielding the matrix values assigned to all possible interactions between pairs of these first five edges.  Using the circumcentric dual mesh and definition (\ref{eqn:hsdiag}), we compute
\begin{equation}
\label{eq:M1diag}
\M_1^{Diag}=\ds\left(\
\begin{array}{cccccc}
\ds \frac{4P^2-1}{4P} & 0 & 0 & 0 & 0 & \cdots \\
\\
0 & \mdiagent & 0 & 0 & 0 & \cdots \\
\\
0 & 0 & \mdiagent & 0 & 0 & \cdots \\
\\
0 & 0 & 0 & \mdiagent & 0 & \cdots \\
\\
0 & 0 & 0 & 0 & \mdiagent & \cdots \\
\vdots & \vdots & \vdots & \vdots & \vdots & \ddots
\end{array}
\right)
\end{equation}
where $\mdiagent=\mdiagentn$.  Since $\M_1^{Diag}$ is diagonal, its condition number is the ratio of its largest diagonal entry to its smallest.  The uncomputed diagonal entries will be very close to $\mdiagent$ meaning the condition number can be approximated as
\[\cond{\M_1^{Diag}}\approx \frac{4P^2-1}{4P}/\mdiagent \in O(P).\]
Using the Whitney interpolant for edges (see (\ref{eq:prwhit1}) in Appendix A) and the definition of $\M_1^{Whit}$ given in (\ref{eqn:hswhit}), we can also compute
\begin{equation}
\label{eq:M1whit}
\M_1^{Whit}=\ds\left(\
\begin{array}{cccccc}
\wabwab & \wabwac & \wabwac & \wabwac & \wabwac & \cdots \\
\\
\wabwac & \wacwac & 0       & \wacwbc & 0       & \cdots \\
\\
\wabwac & 0       & \wacwac & 0       & \wacwbc & \cdots \\
\\
\wabwac & \wacwbc & 0       & \wacwac & 0       & \cdots \\
\\
\wabwac & 0       & \wacwbc & 0       & \wacwac & \cdots \\
\vdots & \vdots & \vdots & \vdots & \vdots & \ddots
\end{array}
\right)
\end{equation}
where $\wabwab=\wabwabn$, $\wabwac=\wabwacn$, $\wacwac=\wacwacn$, and $\wacwbc=\wacwbcn$.  Note that some of the structure of $\M_1^{Whit}$ suggested by (\ref{eq:M1whit}) is an artifice of our ordering of the edges as stated in (\ref{eq:edgeorder}).  However, the remaining diagonal entries of $\M_1^{Whit}$ are all close to $\wacwac$, the entire matrix is symmetric, and the remaining non-zero off-diagonal terms are roughly the same size.  Thus, the eigenvalues of the $5\times 5$ matrix shown in (\ref{eq:M1whit}) allow us to approximate the condition number of $\M_1^{Whit}$.  Using Mathematica, we find analytical expressions for the max and min eigenvalues of the $5\times 5$ matrix and take their ratio to approximate
\[\cond{\M_1^{Whit}}\approx\frac{24P^2+5\sqrt 3 P +\sqrt{288P^4-120\sqrt 3 P^3+3P^2+9}+3}{10\sqrt 3 P + 18}\in O(P)\]

Finally, we compute $(\M_1^{Dual})^{-1}$ using the barycentric dual mesh and definition (\ref{eqn:dhswhit}), yielding
\begin{equation}
\label{eq:M1dual}
\left(\M_1^{Dual}\right)^{-1}=\ds\left(\
\begin{array}{cccccc}
\dwabwab & \dwabwac & \dwabwac & \dwabwac & \dwabwac & \cdots \\
\\
\dwabwac & \dwacwac & \dwacwad & \dwacwbc & 0       & \cdots \\
\\
\dwabwac &  \dwacwad  & \dwacwac & 0   & \dwacwbc & \cdots \\
\\
\dwabwac & \dwacwbc    & 0       & \dwacwac & \dwacwad  & \cdots \\
\\
\dwabwac & 0     & \dwacwbc       & \dwacwad & \dwacwac & \cdots \\
\vdots & \vdots & \vdots & \vdots & \vdots & \ddots
\end{array}
\right)
\end{equation}
where $\dwabwab=\dwabwabn$, $\dwabwac=\dwabwacn$, $\dwacwac=\dwacwacn$, $\dwacwad=\dwacwadn$ and $\dwacwbc=\dwacwbcn$.  Note that analytical expressions of these inner products are not feasible due to the need to compute areas of intersection of irregular polygons in the definition of the $\ol\lambda$ functions.  Instead, using Matlab, we create a simple grid-based quadrature method to estimate the entries of $\left(\M_1^{Dual}\right)^{-1}$ for various values of $P$.  As with $\M_1^{Whit}$, we then estimate the condition number of the entire matrix by the ratio of the max and min eigenvalues of the $5\times 5$ matrix given in (\ref{eq:M1dual}).

The cases $P=2$, 5, and 10 were tested.  The integral required to compute $\dwacwbc$ has support outside of the portion of the dual mesh shown in Figure~\ref{fig:M1compNotn}.  We thus set $\dwacwbc$ to be the same as $\dwabwac$, since both are inner products associated to adjacent edges in the dual mesh.  The computed values of $\dwacwad$ were very small, as expected; we found that setting $\dwacwad$ to zero did not affect the condition number estimate.  Our results are summarized in Table~\ref{tab:condnums}.

\begin{table}[ht]
\centering
\sbox{\strutbox}{\rule{0pt}{0pt}}
\begin{tabular}[.75\textwidth]{@{\extracolsep{\fill}} c|ccc}
%\hline\hline
$P$ & $\cond{\M_1^{Diag}}$ & $\cond{\M_1^{Whit}}$ & $\cond{\left(\M_1^{Dual}\right)^{-1}}$ \\[2mm]
\hline\hline\\[2mm]
2 & 6.3 & 3.2 & 1.5  \\[2mm]
\hline\\[2mm]
5 & 17.2 & 9.9 & 1.3 \\[2mm]
\hline\\[2mm]
10 & 34.6 & 21.6 & 1.4  \\[2mm]
\hline
\end{tabular}
\caption{Comparison of condition numbers of different discrete Hodge stars for various values of $P$.}
\label{tab:condnums}
\end{table}

Our numerical experiments thus provide evidence for the claim
\[\cond{(\M_1^{Dual})^{-1}}\in O(1).\]

The above example confirms that while our dual discrete Hodge star has an analogous definition to the primal discrete Hodge star, its condition number is indeed controlled by the geometric properties of the dual mesh elements, not those of the primal mesh elements.  This fact is especially useful for problems on tetrahedral meshes where slivers (narrow, nearly planar tetrahedra) frequently occur and are difficult to remove.

\section{Applications}
\label{sec:res}

The dual interpolation functions $\ol\I_{n-k}$ we defined in (\ref{eqn:olikdef}) and the dual discrete Hodge star we defined in (\ref{eqn:dhswhit}) are new tools for designing stable finite element methods.  We start by explaining the generic methodology of our approach and then apply it to two sample finite element problems from the literature: magnetostatics and Darcy flow.

\subsection{Generic methodology}

The Discrete Exterior Calculus approach to discretizing a PDE is as follows:
\renewcommand{\labelenumi}{\Roman{enumi}.}
\begin{enumerate}
\item \textbf{Translate} the continuous PDE problem into the language of exterior calculus.
\item \textbf{Linearize} the problem, possibly by introducing an intermediary variable (i.e. a mixed method).
\item \textbf{Discretize} the $k$-forms into $k$-cochains and the operators $d$ and $\ast$ into $\D$ and $\M$ matrices.
\item \textbf{Solve} a linear system constructed from the discrete equations.
\end{enumerate}
Our methodology focuses on step III and exposes how there are often many natural choices for discretization in line with DEC theory.  Consider the case where we are given a PDE in terms of a variable $u$ that is treated as a $k$-form in the continuous setting.   Suppose that a mixed method is possible in which the intermediary variable $v$ should be interpreted as an $n-k-1$ form.  In this case, the typical mixed linear system is
\begin{equation}
\label{eq:gensys1}
\left(\begin{array}{rc}
-\M_k & \D_k^T \\
\D_k & 0
\end{array}
\right)
\left(\begin{array}{c} \pc u \\ \dc v \end{array}
\right)
=
\left(\begin{array}{c} \dc f \\  \pc g \end{array}
\right).
\end{equation}
where $\pc u\in\cs^k$, $\dc v\in\dcs^{n-k-1}$ are the discretized variables and $\dc f\in\dcs^{n-k}$, $\pc g\in\cs^{k+1}$ are the discretized load data.

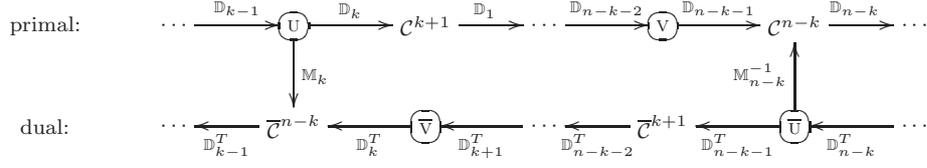
\begin{figure*}
\[
\xymatrix{
\text{primal:}& \cdots \ar[r]^-{\D_{k-1}} & *+[F-:<4pt>]{\pc u}  \ar[d]^{\M_k} \ar[r]^-{\D_{k}} & {\cs^{k+1}} \ar[r]^{\D_1} &
\cdots \ar[r]^-{\D_{n-k-2}} & *+[F-:<4pt>]{\pc v}  \ar[r]^-{\D_{n-k-1}} & {\cs^{n-k}} \ar[r]^{\D_{n-k}} & \cdots \\
\text{dual:} & \cdots & {\dcs^{n-k}}  \ar[l]^-{\D_{k-1}^T} & *+[F-:<4pt>]{\dc v} \ar[l]^-{\D_k^T} & \cdots \ar[l]^-{\D_{k+1}^T} & {\dcs^{k+1}}  \ar[l]^-{\D_{n-k-2}^T} & *+[F-:<4pt>]{\dc u} \ar[u]^{\Mi{n-k}}\ar[l]^-{\D_{n-k-1}^T} & \cdots \ar[l]^-{\D_{n-k}^T}
}
\]
\caption{Portion of a generic DEC-deRham diagram (cf. Figure \ref{fig:derham}) showing the natural duality between the variables and operators of systems (\ref{eq:gensys1}) and (\ref{eq:gensys2}).  Discretizations of the variables are written in place of the primal or dual cochain spaces to which they belong. }
\label{fig:genmethod}
\end{figure*}

The simple idea at the heart of our technique is to swap the \textit{type} of dicretization (primal or dual) of each variable and then infer the rest of the system from DEC theory.  Note that the cochain order of each variable should not change, only the mesh on which it is discretized.  Hence, the dual formulation of system (\ref{eq:gensys1}) is
\begin{equation}
\label{eq:gensys2}
\left(\begin{array}{lc}
-\Mi{n-k} & \D_{n-k-1} \\
\text{ }\D_{n-k-1}^T & 0
\end{array}
\right)
\left(\begin{array}{c} \dc u \\ \pc v \end{array}
\right)
=
\left(\begin{array}{c} \pc f \\  \dc g \end{array}
\right).
\end{equation}
where now $\dc u\in\dcs^k$, $\pc v\in\cs^{n-k-1}$ are the discretized variables and $\pc f\in\cs^{n-k}$, $\dc g\in\dcs^{k+1}$ are the discretized load data.  We show in Figure \ref{fig:genmethod} how these two discretizaions fit into a generic DEC-deRham diagram in a natural and complementary fashion.

Additional equivalent systems can be derived by using proxy variables in clever ways, e.g. solving for some $\pc z\in\cs^{k-1}$ such that $\pc x$ is defined uniquely by $\pc x =\D_{k-1}\pc z$.  These systems are easiest to understand via the specific examples we now examine.

\subsection{Magnetostatics}
\label{subsec:magneto}

The magnetostatics problem is characterized by Gauss's law for magnetism, Amp\`ere's law, and a constitutive relationship, respectively,
\begin{equation}
\label{eq:magcnts}
\Div b=0,\quad  \ast b= h,\quad \Curl h= j.
\end{equation}

Here, $j$ is a given current density and $b$ and $h$ both represent the magnetic field.  It is assumed that the domain $\Omega$ is contractible with boundary $\Gamma$ written as a disjoint union $\Gamma^e\cup \Gamma^h$ such that $\hat n\cdot b=0$ on $\Gamma^e$ and $\hat n\times h=0$ on $\Gamma^h$.

A DEC-based treatment of the problem reveals canonical and symmetrical ways to put this into a mixed formulation linear system, depending on whether $b$ is discretized as a primal or dual cochain.  If we discretize $b$ as a primal 2-cochain $\pc b\in\cs^2$ and $h$ as a dual 1-cochain $\dc h\in\dcs^1$, equations (\ref{eq:magcnts}) become

\[\D_2\pc b = 0,\quad \M_2 \pc b = \dc h,\quad \D_1^T \dc h = \dc j.\]
This allows for two possible mixed systems.  The first is
\begin{equation}
\label{eq:magsys1}
\left(\begin{array}{cc}
-\M_2 & \D_2^T \\
\D_2 & 0
\end{array}
\right)
\left(\begin{array}{c} \pc b \\ \dc p \end{array}
\right)
=
\left(\begin{array}{c} -\dc h_0 \\  0 \end{array}
\right).
\end{equation}
In this system, $\dc h_0\in\dcs^1$ is any dual 1-cochain satisfying  $\D_1^T \dc h_0 = \dc j$ and $\dc h$ is defined by $\dc h:=\dc h_0 + \D_2^T\dc p$.  Thus $\D_1^T\dc h=\D_1^T (\dc h_0 + \D_2^T\dc p) = \dc j$ is assured.

The second mixed system is
\begin{equation}
\label{eq:magsys2}
\left(\begin{array}{cc}
-\Mi 2 & \D_1 \\
\D_1^T & 0
\end{array}
\right)
\left(\begin{array}{c} \dc h \\ \pc a \end{array}
\right)
=
\left(\begin{array}{c} 0 \\ \dc j \end{array}
\right).
\end{equation}
In this system, $\pc b$ is defined by $\pc b := \D_1\pc a$, so that $\D_2\pc b=\D_2\D_1\pc a =0$.  For a fixed $\dc j$, systems (\ref{eq:magsys1}) and (\ref{eq:magsys2}) result in the same solution pair $(\pc b, \dc h)$ and were shown by Bossavit~\cite{B2000part5} to converge to the solution pair $(b,h)$ to (\ref{eq:magcnts}) as the size of mesh elements goes to zero.

We now consider a novel dual discretization approach by treating $b$ as a dual 2-cochain $\dc b\in\dcs^2$ and $h$ as a primal 1-cochain $\pc h\in\cs^1$.  The continuous problem (\ref{eq:magcnts}) is now discretized by
\[\D_0^T\dc b = 0,\quad \dc b = \M_1\pc h,\quad \D_1 \pc h = \pc j.\]
The first mixed system of this dual formulation is
\begin{equation}
\label{eq:magsys3}
\left(\begin{array}{cc}
-\Mi 1 & \D_0 \\
\D_0^T & 0
\end{array}
\right)
\left(\begin{array}{c} \dc b \\ \pc p \end{array}
\right)
=
\left(\begin{array}{c} -\pc h_0 \\ 0 \end{array}
\right).
\end{equation}
In this system, $\pc h_0\in\cs^1$ is any primal 1-cochain satisfying  $\D_1 \pc h_0 = \pc j$ and $\pc h$ is defined by $\pc h:=\Mi 1\dc b$.  Thus $\D_1 \pc h=\D_1 (\pc h_0 + \D_0\pc p) = \pc j$ is assured.  The last system is
\begin{equation}
\label{eq:magsys4}
\left(\begin{array}{cc}
-\M_1 & \D_1^T \\
\D_1 & 0
\end{array}
\right)
\left(\begin{array}{c} \pc h \\ \dc a \end{array}
\right)
=
\left(\begin{array}{c} 0 \\  \pc j \end{array}
\right),
\end{equation}
where $\dc b$ is defined by $\dc b := \D_1^T\dc a$ so that $\D_0^T\dc b = \D_0^T\D_1^T\dc a = 0$.  For a fixed $\pc j$, systems (\ref{eq:magsys3}) and (\ref{eq:magsys4}) will result in the same solution pair $(\dc b, \pc h)$.  In a future work, we will show that these systems also converge to the solution pair $(b,h)$ to (\ref{eq:magcnts}) as the size of mesh elements goes to zero.  Taking that for granted, we state the advantages of having all four systems (\ref{eq:magsys1}), (\ref{eq:magsys2}), (\ref{eq:magsys3}), and (\ref{eq:magsys4}) available for implementation.

First, observe that systems (\ref{eq:magsys1}) and (\ref{eq:magsys2}) make use of the $\M_2$ matrix and its inverse while (\ref{eq:magsys3}) and (\ref{eq:magsys4}) use the $\M_1$ matrix.  If the diagonal Hodge star is used, then $\M_2$ requires good ratios between the size of primal faces and their dual edges while $\M_1$ requires good ratios between the size of primal edges and their dual faces.  Thus, on unstructured meshes, one system may break numerically on a mesh that is acceptable for another system.

Second, if the Whitney Hodge star is used, $\Mi k$ may be a full rank matrix, making systems (\ref{eq:magsys2}) and (\ref{eq:magsys3}) less attractive numerically.  By constructing the dual discrete Hodge stars as proposed in this paper, these systems become sparse again by Lemma~\ref{lem:mkdsparse} and thus are available as a practical alternative.

Third, having four systems available for the same problem allows for rigorous error-checking and cross-confirmation of results.  This is particularly valuable when physical experimental confirmation of the results is impossible or expensive.

\subsection{Darcy Flow}
\label{subsec:dflow}
The Darcy flow problem in $\R^3$ is
\begin{equation}
\label{eq:dfloweqns}
\left\{ \begin{array}{rcll}
 f + \frac k\mu\nabla p & = & 0 & \text{in $\Omega$,} \\
\Div f & = & \phi & \text{in $\Omega$,} \\
f\cdot\hat n & = & \psi & \text{on $\p\Omega$,}
\end{array}\right .
\end{equation}
where $k$ and $\mu$ are physical constants, $f$ is volumetric flux and $p$ is pressure.  It is assumed that there is no external body force, the boundary $\Gamma:=\p\Omega$ is piecewise smooth, and the compatibility condition $\int_\Omega\phi d\Omega = \int_{\p\Omega}\psi d\Gamma$ is satisfied.  Without loss of generality, take $\mu=k$.

First consider discretizing $f$ as a a primal 2-cochain $\pc f\in\cs^2$ and $p$ as a dual 0-cochain $\dc p\in\dcs^0$, yielding the discretized equations
\[\M_2\pc f + \D_2^T \dc p = 0,\quad \D_2\pc f = \Phi.\]
Hirani el al.~\cite{HNC2008} used this approach to derive the linear system
\begin{equation}
\label{eqn:dfsys1}
\left(\begin{array}{cc}
\M_2 & \D_2^T \\
\D_2 & 0
\end{array}
\right)
\left(\begin{array}{c} \pc f \\ \dc p \end{array}
\right)
=
\left(\begin{array}{c} 0 \\ \Phi\end{array}
\right).
\end{equation}
We present an alternative formulation using the same discretization, inspired by the magnetostatics systems (\ref{eq:magsys2}) and (\ref{eq:magsys4}).  Let $\pc f_0\in \cs^2$ be a primal 2-cochain satisfying $\D_2\pc f_0=\Phi$.  The system is

\begin{equation}
\label{eqn:dfsys2}
\left(\begin{array}{cc}
-\Mi 2 & \D_1 \\
\D_1^T & 0
\end{array}
\right)
\left(\begin{array}{c} \dc q \\ \pc g \end{array}
\right)
=
\left(\begin{array}{c} -\pc f_0 \\ 0\end{array}
\right).
\end{equation}
Here, $\dc p$ is a solution to $\D_2^T\dc p=-\dc q$.  The existence of $\dc p$ is guaranteed by the exactness of the dual cochain sequence at $\dcs^1$ and uniqueness of $\dc p$ is determined by initial conditions or boundary data.  The flux cochain $\pc f$ is defined to be $\Mi 2\dc q$ so that $\D_2\pc f = \D_2(\Mi 2\dc q)=\D_2(\pc f_0+\D_1\pc g)=\phi$.

We now present the dual formulations derived by treating $f$ as a dual 2-cochain $\dc f\in\dcs^2$ and $p$ as a primal 0-cochain $\pc p\in\cs^0$.  The discretized equations are now
\[\Mi 1\dc f + \D_0 \pc p = 0,\quad \D_0^T\dc f = \ol\Phi.\]
The first system of this formulation is
\begin{equation}
\label{eqn:dfsys3}
\left(\begin{array}{cc}
\Mi 1 & \D_0 \\
\D_0^T & 0
\end{array}
\right)
\left(\begin{array}{c} \dc f \\ \pc p \end{array}
\right)
=
\left(\begin{array}{c} 0 \\ \ol\Phi\end{array}
\right).
\end{equation}
The second system is
\begin{equation}
\label{eqn:dfsys4}
\left(\begin{array}{cc}
\M_1 & \D_1^T \\
\D_1 & 0
\end{array}
\right)
\left(\begin{array}{c} \pc q \\ \dc g \end{array}
\right)
=
\left(\begin{array}{c} \dc f_0 \\ 0\end{array}
\right).
\end{equation}
where $\dc f_0$ is a solution to $\D_0^T\dc f_0=\ol\Phi$ and $\dc f$ is defined to be $\M_1 \pc q$, analogous to system (\ref{eqn:dfsys2}).  Thus, taking $\D_0^T$ of both sides of the top equation of (\ref{eqn:dfsys4}) yields $\D_0^T\dc f = \ol\Phi$.  Further, the bottom equation of (\ref{eqn:dfsys4}) yields $\D_1\pc q = 0$ which, by the exactness property of the primal cochain sequence implies that there exists a solution $\pc p$ to $\D_0\pc p=-\pc q$.

We now have four mixed systems, (\ref{eqn:dfsys1})-(\ref{eqn:dfsys4}), discretizing the Darcy flow equations (\ref{eq:dfloweqns}), three of which had not be considered by Hirani et al.~\cite{HNC2008}.  This plethora of equivalent systems offers the same advantages as those discussed at the end of the magnetostatics example from Section~\ref{subsec:magneto}.

\section{Conclusion}

In this work we have augmented the theories of Discrete Exterior Calculus and mixed methods by introducing two novel tools: Whitney-like interpolation functions defined on dual domain meshes and a sparse inverse discrete Hodge star.  We have shown the tools to have natural, straightforward definitions and clear geometric interpretations.  We have used them to derive previously unexamined numerical stability criteria relating to the condition number of the discrete Hodge star used in the method, based on the geometry of the dual mesh cells.  Further, we have demonstrated in both general and specific contexts how these tools can be used to develop alternative discretizations of PDEs with sparse, well-conditioned matrices.  The techniques we have described provide a valuable methodology for researchers to revisit their current finite element formulations and confirm or improve their results with new discretization methods.

\appendix

\section{Whitney Functions for Primal Meshes}

Whitney $k$-forms are piecewise linear functions on a primal mesh, one for each $k$-simplex in the mesh.
\begin{itemize}
\item\textbf{Primal Vertices.} The Whitney 0-form associated to a vertex $\sigma^0:= \bv_i$ is denoted
\[\whit_{\sigma^0}:=\lambda_i,\]
where $\lambda_i$ is the barycentric function for the vertex.  More precisely, $\lambda_i$ is defined by the condition of being linear on every simplex of the mesh, subject to the constraints $\lambda_i(\bv_j)=\delta_{ij}$.
\item\textbf{Primal Edges.} The Whitney 1-form associated to an oriented edge $\sigma^1:=[\bv_i,\bv_j]$ is the vector-valued function
\begin{equation}
\label{eq:prwhit1}
\whit_{\sigma^1}:=\lambda_i\nabla\lambda_j-\lambda_j\nabla\lambda_i.
\end{equation}
\item\textbf{Primal Faces.}
The Whitney 2-form associated to an oriented face $\sigma^2:=[\bv_i,\bv_j,\bv_k]$ is the vector-valued function
\begin{equation}
\label{eq:prwhit2}
\whit_{\sigma^2}:= 2\left(\whitC i j k+\whitC j ki + \whitC k i j\right)
\end{equation}
\item\textbf{Primal Tetrahedra.}\footnote{Note that the $\whit_{\sigma^3}$ definition has been simplified from a more general definition of Whitney forms~\cite{W1957} using the geometric identity
\[\nabla\lambda_i\cdot(\nabla\lambda_j\times\nabla\lambda_k)=\pm \frac 1{3!|\sigma^3|}\]
where the right side has sign $-1$ if an odd index was omitted from the scalar triple product and $+1$ otherwise. This reduces the sum in the general formula to $(1/|\sigma^3|)\sum_i\lambda_i$, which is simply  $1/|\sigma^3|$ due to the partition of unity formed by the barycentric functions.}
The Whitney 3-form associated to an oriented tetrahedron $\sigma^3$ is its characteristic function, scaled by the reciprocal of the volume $\sigma^3$.
\[\whit_{\sigma^3}:= \chi_{\sigma^3} = \left\{\begin{array}{rl} 1/|\sigma^3| & \text{on $\sigma^3$}\\ 0 & \text{otherwise}  \end{array}\right.  \]

\end{itemize}

\section{Generalized Barycentric Functions}

Let $\Tau$ be a top-dimensional cell of the dual mesh (i.e. a polygon in 2D or a polyhedron in 3D) with vertices $\bv_1,\ldots,\bv_N$.     A set of functions $\ol\lambda_i:\Tau\raw\R$, $i=1,\ldots, N$ are called \textbf{barycentric coordinates} on $\Tau$ if they satisfy two properties.
\renewcommand{\labelenumi}{B\arabic{enumi}.}

\begin{enumerate}
\item \textbf{Non-negative}:  $\ol\lambda_i\geq 0$.
\item \textbf{Linear Completeness}: For any linear function $L:\Tau\raw\R$,
\[\ds L=\sum_{i=1}^{N} L(\bv_i)\ol\lambda_i.\]
\end{enumerate}
A set of barycentric coordinates $\{\ol\lambda_i\}$ also satisfies these additional familiar properties:
\renewcommand{\labelenumi}{B\arabic{enumi}.}
\begin{enumerate}
\setcounter{enumi}{2}
\item \textbf{Partition of unity:} $\ds\sum_{i=1}^{N}\ol\lambda_i\equiv 1$. \label{b:partition}
\item \textbf{Linear precision:} $\ds\sum_{i=1}^{N}\bv_i\ol\lambda_i(\bx)=\bx$. \label{b:linprec}
\item \textbf{Interpolation:} $\ds\ol\lambda_i(\bv_j) = \delta_{ij}$. \label{b:interpolation}
\end{enumerate}
A proof that properties B3-B5 are implied by B1-B2 in the 2D case can be found in our paper~\cite{GRB2010}.  The 3D case is similar.

Three major approaches to defining generalized barycentric functions on 2D polygons have emerged in the literature.  The Wachspress functions~\cite{W1975,FHK2006} are rational functions constructed explicitly based on the areas of certain triangles within $\Tau$.  The Sibson functions~\cite{S1980}, also called the natural neighbor or natural element coordinates~\cite{SM2006}, are also constructed explicitly, but instead use the areas of Voronoi regions associated with the vertices of $\Tau$.  The Harmonic functions~\cite{WSHD2007,C2008} are defined as the solution to Laplace's equation over $\Tau$ with certain piecewise linear boundary data.

We have shown in~\cite{GRB2010} that any of these functions suffice to give the optimal interpolation estimate for the lowest order case in 2D, assuming some basic geometric quality criteria on the dual mesh elements.  For this paper, we have employed only the Sibson coordinates as they generalize easily to 3D, are reasonable to implement, and are more stable against bad geometry than the Wachspress functions.  A formal proof of their convergence properties in 3D will be the focus of a future work.


%% References

\paragraph{Acknowledgments} We are grateful to Alexander Rand for his help in implementing the Sibson coordinates.  This research was supported in part by NIH contracts R01-EB00487, R01-GM074258, and a grant from the UT-Portugal CoLab project.


%\bibliographystyle{model3-num-names}
\bibliographystyle{abbrv}
\bibliography{dfmfem-arxiv}

\end{document}